%% file: 2003-16.tex
\documentclass{gtart}
\input gtspec

\lognumber{267}
\volumenumber{7}\papernumber{16}
\volumeyear{2003}
\pagenumbers{569}{613}
\received{2 August 2002}
\accepted{9 October 2003}
\published{11 October 2003}

\proposed{Steve Ferry}

\seconded{Martin Bridson, Frances Kirwan}

\usepackage{amsmath,amsfonts,amscd,amssymb}

\DeclareMathAlphabet\EuR{U}{eur}{m}{n}
\SetMathAlphabet\EuR{bold}{U}{eur}{b}{n}
\usepackage[dvips]{graphicx} 
\usepackage[arrow,matrix]{xy}

\newtheorem{thm}{Theorem}[section]   
\newtheorem{lem}[thm]{Lemma}         
\newtheorem{prop}[thm]{Proposition}

\theoremstyle{definition}
\newtheorem{defn}[thm]{Definition}  
\theoremstyle{remark}
\newtheorem{rem}[thm]{Remark}

\newtheorem{notation}[thm]{Notation}
{\catcode`@=11\global\let\c@equation=\c@thm}

\renewcommand{\theenumi}{\alph{enumi}}

\newcommand{\squarematrix}[4]{\left( \begin{array}{cc} #1 & #2 \\ #3 &
#4
\end{array} \right)}

\newcommand{\comsquare}[8]{
\begin{center}
$\begin{CD}
#1 @>#2>> #3\\
@V{#4}VV @VV{#5}V\\
#6 @>>#7> #8
\end{CD}$
\end{center}}

\newcommand{\curs}{\EuR}
\newcommand{\Chain}{\curs{Chain}}

\newcommand{\Mod}{\curs{Mod}}
\newcommand{\Or}{\curs{Or}}

\newcommand{\Orcat}{\Or(G;\calfin)}

\newcommand\bZ{\mathbb Z}
\newcommand\bR{\mathbb R}

\newcommand{\zz}{{\mathbb Z}}

\newcommand{\cc}{{\mathbb C}}
\newcommand{\hh}{{\mathbb H}}

\newcommand{\qq}{{\mathbb Q}}

\newcommand{\rr}{{\mathbb R}}

\newcommand\cH{\mathcal H}
\newcommand\cE{\mathcal E}
\newcommand\cL{\mathcal L}
\newcommand\cK{\mathcal K}

\newcommand{\calfin}{\mathcal{F}\curs{in}}

\newcommand{\caln}{{\cal N}}

\newcommand\colim{\raisebox{-1ex}{\ensuremath{
\stackrel{\textstyle{\lim}}{\scriptstyle{\longrightarrow}}}}}
\newcommand\Cl{\operatorname{Cliff}(TM)}
\newcommand\dvol{d\text{vol}}
\newcommand\range{\operatorname{range}}
\newcommand\dom{\operatorname{dom}}
\newcommand\lp{\textup{(}}
\newcommand\rp{\textup{)}}

\newcommand{\asmb}{\operatorname{asmb}}
\newcommand{\aut}{\operatorname{aut}}
\newcommand{\ch}{\operatorname{ch}}
\newcommand{\chern}{\operatorname{chern}}

\newcommand{\class}{\operatorname{class}}

\newcommand{\con}{\operatorname{con}}

\newcommand{\consub}{\operatorname{consub}}

\newcommand{\Eul}{\operatorname{Eul}}

\newcommand{\Gen}{\operatorname{Gen}}
\newcommand{\HS}{\operatorname{HS}}

\newcommand{\id}{\operatorname{id}}
\newcommand{\ind}{\operatorname{ind}}

\newcommand{\Is}{\operatorname{Is}}

\newcommand{\map}{\operatorname{map}}
\newcommand{\MOD}{\Mod}

\newcommand{\mor}{\operatorname{mor}}

\newcommand{\pr}{\operatorname{pr}}
\newcommand{\pt}{\text{\textup{pt}}}
\newcommand{\Rep}{\operatorname{Rep}}

\newcommand{\rank}{\operatorname{rank}}

\newcommand{\Sign}{\operatorname{Sign}}

\newcommand{\supp}{\operatorname{supp}}

\newcommand{\tr}{\operatorname{tr}}

\newcommand{\Zero}{\operatorname{Zero}}

\newcommand{\tit}[2]{\section{#1}
\label{sec: #2}}
\newcommand{\subtit}[2]{\subsection{#1}
\label{subsec: #2}}

\newenvironment{Relax}{\relax}{\relax}

\begin{document}
\title[Equivariant Euler characteristics and $K$-homology]{Equivariant
Euler characteristics and $K$-homology\\ Euler classes
for proper cocompact $G$-manifolds}                    
\asciititle{Equivariant
Euler characteristics and K-homology Euler classes
for proper cocompact G-manifolds}                    

\authors{Wolfgang L\"uck\\Jonathan Rosenberg} 
\asciiauthors{Wolfgang Lueck and Jonathan Rosenberg}
\coverauthors{Wolfgang L\noexpand\"uck\\Jonathan Rosenberg}
\address{Institut f\"ur Mathematik und Informatik,
Westf\"alische Wilhelms-Universtit\"at\\
Einsteinstr.\ 62,
48149 M\"unster, Germany}
\secondaddress{Department of Mathematics, University of 
Maryland\\College
Park, MD 20742, USA}
\asciiaddress{Institut fur Mathematik und Informatik,
Westfalische Wilhelms-Universtitat\\Einsteinstr. 62,
48149 Munster, Germany\\and\\Department of Mathematics, 
University of 
Maryland\\College
Park, MD 20742, USA}

\email{lueck@math.uni-muenster.de}
\secondemail{jmr@math.umd.edu}
\asciiemail{lueck@math.uni-muenster.de, jmr@math.umd.edu}
\urladdr{wwwmath.uni-muenster.de/u/lueck,
www.math.umd.edu/\raisebox{-.6ex}{\symbol{"7E}}jmr}


\begin{abstract}   
Let $G$ be a countable discrete group and let $M$ be a smooth proper
cocompact $G$-manifold without boundary. The Euler operator defines via
Kasparov theory an element, called the equivariant Euler class,  
in the equivariant $KO$-homology of $M\!$. The universal equivariant Euler
characteristic of $M\!$, which lives in a group $U^G(M)$,
counts the equivariant cells of $M$, taking the
component structure of the various fixed point sets into account.
We construct a natural homomorphism from $U^G(M)$ to the
equivariant $KO$-homology of $M\!$. The main result of this paper says
that this map sends the universal equivariant Euler characteristic to
the equivariant Euler class. In particular this shows that there are
no ``higher'' equivariant Euler characteristics. We show that,
rationally, the equivariant Euler class 
carries the same information as the collection of the orbifold Euler 
characteristics of the components of the $L$-fixed point sets $M^L\!$,
where $L$ runs through the finite cyclic subgroups of $G$.  
However, we give an example of an action of the symmetric
group $S_3$ on the $3$-sphere for which the equivariant Euler class
has order $2$, so there is also some torsion information.
\end{abstract}

\asciiabstract{Let G be a countable discrete group and let M be a
smooth proper cocompact G-manifold without boundary. The Euler
operator defines via Kasparov theory an element, called the
equivariant Euler class, in the equivariant K-homology of M. The
universal equivariant Euler characteristic of M, which lives in a
group U^G(M), counts the equivariant cells of M, taking the component
structure of the various fixed point sets into account.  We construct
a natural homomorphism from U^G(M) to the equivariant KO-homology of
M. The main result of this paper says that this map sends the
universal equivariant Euler characteristic to the equivariant Euler
class. In particular this shows that there are no `higher' equivariant
Euler characteristics. We show that, rationally, the equivariant Euler
class carries the same information as the collection of the orbifold
Euler characteristics of the components of the L-fixed point sets M^L,
where L runs through the finite cyclic subgroups of G.  However, we
give an example of an action of the symmetric group S_3 on the
3-sphere for which the equivariant Euler class has order 2, so there
is also some torsion information.}

\begin{Relax}\end{Relax}

              \


\asciikeywords{Equivariant K-homology, de Rham operator, signature
operator, Kasparov theory, equivariant Euler characteristic, fixed
sets, cyclic subgroups, Burnside ring, Euler operator, equivariant
Euler class, universal equivariant Euler characteristic}

{\small\maketitlepage}
%

\setcounter{section}{-1}
\tit{Background and statements of results}{Background and statements
of results} 

Given a countable discrete group $G$ and a cocompact proper smooth
$G$-manifold $M$ without boundary and with $G$-invariant Riemannian
metric, the Euler characteristic operator
defines via Kasparov theory
an element, the \emph{equivariant Euler class},  in the equivariant real $K$-homology group of $M$
\begin{eqnarray}
\Eul^G(M) & \in & KO_0^G(M). \label{definition of equivariant Euler class}
\end{eqnarray}
The \emph{Euler characteristic operator} is the minimal closure, or
equivalently, the maximal closure, of the densely defined operator
$$(d+d^*)\co \Omega^*(M) \subseteq L^2\Omega^*(M) \to L^2\Omega^*(M),$$
with the $\zz/2$-grading coming from the degree of a differential
$p$-form. The \emph{equivariant signature operator},
defined when the manifold is equipped with a $G$-invariant
orientation, is the same
underlying operator, but with a different grading coming from the 
Hodge star operator. The signature operator also defines an element
\begin{eqnarray*}
\Sign^G(M) & \in & K_0^G(M),
\end{eqnarray*}
which carries a lot of geometric information about the action of
$G$ on $M$. (Rationally, when $G=\{1\}$, $\Sign(M)$
is the Poincar\'e dual of the total $\mathcal L$-class,
the Atiyah-Singer $L$-class, which differs from the Hirzebruch
$L$-class only by certain well-understood powers of $2$,
but in addition, it also carries quite interesting integral
information \cite{KM}, \cite{Ro2}, \cite{RW2}. A partial
analysis of the class $\Sign^G(M)$ for $G$ finite may be found in
\cite{RW} and \cite{Rosenberg(1999b)}.)

We want to study how much information $\Eul^G(M)$ carries.
This has already been done by the second author 
\cite{Rosenberg(1999a)} in the non-equivariant case. Namely, 
given a closed  Riemannian manifold $M$, not necessarily connected, let
$$e\co \bigoplus_{\pi_0(M)} \zz = \bigoplus_{\pi_0(M)}
KO_0(\{\ast\}) ~ \to ~ KO_0(M)$$ 
be the map induced by the various inclusions $\{\ast\} \to M$. This map is
split injective; a splitting is given by the various projections $C \to
\{\ast\}$ for $C \in \pi_0(M)$, and  sends $\{\chi(C) \mid C \in
\pi_0(M)\}$ to $\Eul(M)$. 
Hence $\Eul(M)$ carries precisely the same information as the Euler
characteristics of the various components of $M\!$, 
and there are no ``higher'' Euler classes. Thus
the situation is totally different from what happens with the signature
operator.

We will see that in the equivariant case there are again no ``higher''
Euler characteristics and that $\Eul^G(M)$ is determined by the 
\emph{universal equivariant Euler characteristic} (see Definition 
\ref{def: universal equivariant Euler characteristic})
\begin{eqnarray*} 
\chi^G(M) & \in & U^G(M) = \bigoplus_{(H) \in \consub(G)}~
\bigoplus_{W\!H\backslash \pi_0(M^H)} \zz.
\end{eqnarray*}
Here and elsewhere $\consub(G)$ is the set of conjugacy classes of subgroups
of $G$ and $N\!H = \{g \in G \mid g^{-1}Hg = H\}$ is the
\emph{normalizer} of the subgroup 
$H\subseteq G$ and $W\!H:= N\!H/H$ is its \emph{Weyl group}.
The component of $\chi^G(M)$ associated to $(H) \in
\consub(G)$ and $W\!H\cdot C \in W\!H\backslash \pi_0(M^H)$ is the
(ordinary) Euler characteristic 
$\chi(W\!H_C\backslash (C,C\cap M^{>H}))$, where $W\!H_C$ is the isotropy
group of $C \in \pi_0(M^H)$ under the $W\!H$-action.
There is a natural homomorphism
\begin{eqnarray}
e^G(M)\co U^G(M) & \to & KO_0^G(M). \label{definition of e^G(M)}
\end{eqnarray}
It sends the basis element associated to  $(H) \subseteq \consub(G)$ 
and $W\!H \cdot C \in W\!H\backslash \pi_0(M^H)$ to the image
of the class of the trivial $H$-representation $\rr$ under the composition
$$R_{\rr}(H) = KO_0^H(\{\ast\}) \xrightarrow{(\alpha)_*} KO_0^G(G/H)
\xrightarrow{KO_0^G(x)} KO_0^G(M),$$  
where $(\alpha)_*$ is the isomorphism coming from induction via the
inclusion $\alpha\co H \linebreak\to G$ and 
$x\co G/H \to M$ is any $G$-map with $x(1H) \in C$. 
The main result of this paper is

\begin{thm}[Equivariant Euler class and Euler characteristic]
\label{the: e^G(M)(chi^G(M)) = Eul^G(M)}
Let $G$ be a countable discrete group and let 
$M$ be a cocompact proper smooth $G$-manifold without boundary. Then
$$e^G(M)(\chi^G(M)) = \Eul^G(M).$$
\end{thm}

The proof of Theorem \ref{the: e^G(M)(chi^G(M)) = Eul^G(M)} involves
two independent steps. Let $\Xi$ be an equivariant vector field on $M$
which is 
transverse to the zero-section. Let $\Zero(\Xi)$ be the set
of points $x \in M$ with $\Xi(x) = 0$. Then $G\backslash \Zero(\Xi)$ is
finite. The zero-section $i\co M \to TM$
and the inclusion  $j_x\co T_xM \to TM$ induce an isomorphism of
$G_x$-representations 
$$T_xi \oplus T_0j_x\co T_xM \oplus T_xM \xrightarrow{\cong}  T_{i(x)}(TM)$$
if we identify $T_0(T_xM) = T_xM$ in the obvious way.  If $\pr_i$
denotes the projection onto the $i$-th factor for $i = 1,2$ we obtain a
linear $G_x$-equivariant isomorphism
\begin{eqnarray}
\hspace{-5mm} &d_x\Xi\co 
T_xM ~ \xrightarrow{T_x\Xi} T_{i(x)}(TM) ~
~ \xrightarrow{(T_xi \oplus T_xj_x)^{-1}} ~
T_xM \oplus T_xM ~\xrightarrow{\pr_2} ~ T_xM. &
\label{definition of d_xv}
\end{eqnarray} 
Notice that we obtain the identity if we replace $\pr_2$ by $\pr_1$ in
the expression \eqref{definition of d_xv} above. 
One can even achieve that $\Xi$ is \emph{canonically transverse to the
zero-section}, i.e., it is transverse to the zero-section 
and $d_x\Xi$ induces the identity on $T_xM/(T_xM)^{G_x}$ for
$G_x$ the isotropy group of $x$ under the $G$-action. This is proved in
\cite[Theorem 1A on page 133]{Waner-Wu(1986)} in the case of a finite group
and the argument directly carries over to the proper cocompact case.
Define the \emph{index} of $\Xi$ at a zero $x$ by
$$s(\Xi,x) ~ = ~ \frac{\det\left((d_x\Xi )^{G_x}\co (T_xM)^{G_x} \to
(T_xM)^{G_x}\right)} 
{\left|\det\left((d_x\Xi )^{G_x}\co (T_xM)^{G_x} \to
(T_xM)^{G_x}\right)\right|} 
\hspace{5mm} \in \{\pm 1\}.$$
For $x \in M$ let $\alpha_x\co G_x \to G$ be the inclusion,
$(\alpha_x)_* \co R_{\rr}(G_x) = KO_0^{G_x}(\{\ast\}) \to
KO_0^G(G/G_x)$ be the map 
induced by induction via  $\alpha_x$ and let $x\co G/G_x \to M$ be the
$G$-map sending $g$ to $g\cdot x$. 
By perturbing the equivariant Euler operator using the vector field
$\Xi$ we will show :

\begin{thm} [Equivariant Euler class and vector fields]
\label{the: Eul^G(X) and vector fields}
Let $G$ be a countable discrete group and let 
$M$ be a cocompact proper smooth $G$-manifold without boundary.
Let $\Xi$ be an equivariant vector field which is canonically
transverse to the zero-section. Then
$$\Eul^G(M) ~ = ~ \sum_{Gx \in G\backslash \Zero(\Xi)}~ 
s(\Xi,x) \cdot KO_0^G(x) \circ (\alpha_x)_*([\rr]),$$
where $[\rr] \in R_{\rr}(G_x) = K_0^{G_x}(\{\ast\})$ is the class of
the trivial $G_x$-representation 
$\rr$, we consider $x$ as a $G$-map $G/G_x \to M$ and $\alpha_x\co
G_x \to G$ is the inclusion. 
\end{thm}

In the second step one has to prove
\begin{eqnarray}
e^G(M)(\chi^G(M)) & = & \sum_{Gx \in G\backslash \Zero(\Xi)}~ 
s(\Xi,x) \cdot  KO_0^G(x) \circ (\alpha_x)_*([\rr]). 
\label{e^G(M)(chi^G(M)}
\end{eqnarray}
This is a direct conclusion of the equivariant Poincar\'e-Hopf theorem
proved in \cite[Theorem 6.6]{Lueck-Rosenberg(2002a)} (in turn a
consequence of the equivariant Lefschetz fixed point
theorem proved in \cite[Theorem 0.2]{Lueck-Rosenberg(2002a)}), which says 
\begin{eqnarray}
\chi^G(M) & = & i^G(\Xi). 
\label{chi^G(M) = i^G(Xi)}
\end{eqnarray}
where $i^G(\Xi)$ is the equivariant index of the vector field $\Xi$
defined in \cite[(6.5)]{Lueck-Rosenberg(2002a)}.
Since we get directly from the definitions
\begin{eqnarray}
e^G(M)(i^G(\Xi)) & = & \sum_{Gx \in G\backslash \Zero(\Xi)}~ 
s(\Xi,x) \cdot  KO_0^G(x) \circ (\alpha_x)_*([\rr]),
\label{e^G(M)(i^G(v))}
\end{eqnarray}
equation \eqref{e^G(M)(chi^G(M)} follows from 
\eqref{chi^G(M) = i^G(Xi)} and  \eqref{e^G(M)(i^G(v))}. 
Hence Theorem
\ref{the: e^G(M)(chi^G(M)) = Eul^G(M)} is true if we can prove
Theorem \ref{the: Eul^G(X) and vector fields}, which will be done in
Section \ref{sec: Perturbing the equivariant Euler operator by a
vector field}. 

We will factorize $e^G(M)$ as
\begin{multline*}
e^G(M)\co U^G(M) \xrightarrow{e_1^G(M)} H_0^{\Or(G)}(M;\underline{R_{\qq}})
\xrightarrow{e_2^G(M)} H_0^{\Or(G)}(M;\underline{R_{\rr}})
\\
\xrightarrow{e_3^G(M)} KO_0^G(M),
\end{multline*}
where $H_0^{\Or(G)}(M;\underline{R_{F}})$ is the Bredon homology of $M$ with 
coefficients in the coefficient system which sends $G/H$ to the 
representation ring 
$R_{F}(H)$ for the field $F = \qq, \rr$. We will show that $e_2^G(M)$
and $e_3^G(M)$  
are rationally injective (see Theorem \ref{the: e^G(X)}).  
We will analyze the map
$e_1^G(M)$, which is not rationally injective in general, in
Theorem \ref{the: commutative diagram for e^G_1(M)}.

The rational information carried by $\Eul^G(M)$ can be expressed in terms of
orbifold Euler characteristics of the various components of the
$L$-fixed point sets 
for all finite cyclic subgroups $L \subseteq G$.  For a component
$C \in \pi_0(M^H)$ denote by $W\!H_C$ its isotropy group under the
$W\!H$-action on 
$\pi_0(M^H)$. For $H \subseteq G$ finite $W\!H_C$ acts properly and
cocompactly on $C$ and its  
\emph{orbifold Euler characteristic} 
(see Definition \ref{def: universal equivariant Euler characteristic}),
which agrees with the more general notion of \emph{$L^2$-Euler
characteristic}, 
$$\chi^{\qq W\!H_C}(C) \in \qq,$$
is defined. Notice that for finite $W\!H_C$ the orbifold Euler
characteristic is given in  
terms of the ordinary Euler characteristic by
$$\chi^{\qq W\!H_C}(C)  = \frac{\chi(C)}{|W\!H_C|}.$$
There is a character map (see \eqref{definition of character map chi^G(X)})
$$\ch^G(M)\co U^G(M) \to \bigoplus_{(H) \in \consub(G)} ~ 
\bigoplus_{W\!H \backslash \pi_0(M^H)} \qq$$
which sends $\chi^G(M)$ to the various $L^2$-Euler 
characteristics $\chi^{\qq W\!H_C}(C)$ for
$(H) \in \consub(G)$ and $W\!H \cdot C \in W\!H\backslash \pi_0(M^L)$. 
Recall that rationally
$\Eul^G(M)$ carries the same information as $e^G_1(M)(\chi^G(M))$
since the rationally injective map $e^G_3(M) \circ e^G_2(M)$ sends 
$e^G_1(M)(\chi^G(M))$ to $\Eul^G(M)$. Rationally 
$e^G_1(M)(\chi^G(M))$ is the same as the collection of all these
orbifold Euler characteristics $\chi^{\qq W\!H_C}(C) $
if one restricts to finite cyclic subgroups $H$. Namely, we will prove
(see Theorem \ref{the: commutative diagram for e^G_1(M)}):

\begin{thm} \label{the: rational information carried by Eul^G(M)}
There is a bijective natural map
$$\gamma^G_{\qq}\co 
\bigoplus_{\substack{(L) \in \consub(G)\\L \text{\textup{ finite cyclic}}}}~  
\bigoplus_{W\!L\backslash \pi_0(M^L)} ~ \qq
~ \xrightarrow{\cong} ~ \qq \otimes_{\zz}
H_0^{\Or(G)}(X;\underline{R_{\qq}})$$ 
which maps 
$$\{\chi^{\qq W\!L_C)}(C) \mid (L) \in \consub(G), L \text{\textup{
finite cyclic}},  
W\!L \cdot C \in W\!L\backslash \pi_0(M^L)\}$$
to $1 \otimes_{\zz} e^G_1(M)(\chi^G(M))$. 
\end{thm}

However, we will show that $\Eul^G(M)$ does carry some torsion information.
Namely, we will prove:
\begin{thm} \label{the: example of Eul^G(M) of order two} 
There exists an action of the symmetric group $S_3$ of order $3!$ on
the $3$-sphere $S^3$ such that 
$\Eul^{S_3}(S^3) \in KO_0^{S_3}(S^3)$ has order $2$.
\end{thm}

The relationship between $\Eul^G(M)$ and the various notions of
equivariant Euler characteristic is clarified in sections
\ref{sec: Review of notions of equivariant Euler characteristic} and
\ref{subsec: Independence of stable equivariant Euler 
characteristic}. 

The paper is organized as follows:

\smallskip
\begin{tabular}{ll}
\ref{sec: Perturbing the equivariant Euler operator by a vector field}
& Perturbing the equivariant Euler operator by a vector field
\\
\ref{sec: Review of notions of equivariant Euler characteristic} & 
Review of notions of equivariant Euler characteristic
\\
\ref{sec: The transformation eG(M)} & The transformation $e^G(M)$
\\
\ref{sec: Examples} & Examples
\\
& \ref{subsec: Finite groups and connected non-empty fixed point sets}
~ Finite groups and connected non-empty fixed point sets
\\
& \ref{subsec: EulG(M) carries torsion information}
~
The equivariant Euler class carries torsion information
\\
& \ref{subsec: Independence of stable equivariant Euler
characteristic} 
~
Independence of $\Eul^G(M)$ and $\chi^G_s(M)$
\\
& \ref{subsec: The image of EulG(M) under the assembly maps}
~ The image of the equivariant Euler class under assembly\\
 & References
\end{tabular}
\smallskip

This paper subsumes and replaces the preprint \cite{Rosenberg(unpub)},
which gave a much weaker version of Theorem 
\ref{the: e^G(M)(chi^G(M)) = Eul^G(M)}.

This research was supported by Sonderforschungsbereich 478
(``Geometrische Strukturen in der Mathematik'') of the
University of M\"unster.
Jonathan Rosenberg was partially supported by NSF grants 
DMS-9625336 and DMS-0103647.

\tit{Perturbing the equivariant Euler operator by a vector field}
{Perturbing the equivariant Euler operator by a vector field}

Let $M^n$ be a complete Riemannian manifold without
boundary, equipped with an isometric action of a discrete group $G$. 
Recall that the de~Rham operator $D=d+d^*$, acting on differential
forms on $M$ (of all possible degrees) is a formally
self-adjoint elliptic operator, and that on the
Hilbert space of $L^2$ forms, it is essentially
self-adjoint \cite{Ga}. With a certain grading on the form bundle
(coming from the Hodge $*$-operator), $D$ becomes the
\emph{signature operator}; with the more obvious grading
of forms by parity of the degree, $D$ becomes the
\emph{Euler characteristic operator} or simply the
\emph{Euler operator}. When $M$ is compact and $G$ is finite,
the kernel of $D$, the space of harmonic forms, is naturally identified 
with the real or complex\footnote{depending on what scalars one is using}
cohomology of $M$ by
the Hodge Theorem, and in this way one observes that
the (equivariant) index of $D$ (with respect to the parity grading) 
in the real representation ring of $G$ is
simply the (equivariant homological) Euler characteristic of $M$, whereas the
index with respect to the other grading is the $G$-signature
\cite{AS3}.

Now by Kasparov theory (good general references are
\cite{Bl} and \cite{Hig}; for the detailed original papers, see
\cite{K} and \cite{K1}), an elliptic operator such as
$D$ gives rise to an equivariant $K$-homology class. In the case of
a compact manifold, the equivariant index of the operator is recovered
by looking at the image of this class under the map collapsing $M$ 
to a point. However, the $K$-homology
class usually carries far more information than the index
alone; for example, it determines the $G$-index of the
operator with coefficients in any $G$-vector bundle, and even determines the
families index in $K^*_G(Y)$ of a family of twists of the operator, 
as determined by a $G$-vector bundle on $M\times Y$. ($Y$ here is an
auxiliary parameter space.) When $M$ is non-compact, 
things are similar, except that usually there is no index, and
the class lives in an appropriate Kasparov
group $K^{-*}_G(C_0(M))$, which is \emph{locally finite} $K^G$-homology,
\emph{i.e.}, the relative group $K_*^G(\overline M, \{\infty\})$, 
where $\overline M$ is the one-point compactification of $M$.\footnote{Here
$C_0(M)$ denotes continuous real- or complex-valued functions on $M$
vanishing at infinity, depending on whether one is using real or complex
scalars. This algebra is \emph{contra}variant in $M$, so a contravariant
functor of $C_0(M)$ is \emph{co}variant in $M$.
Excision in Kasparov theory identifies $K^{-*}_G(C_0(M))$
with $K^{-*}_G(C(\overline M),\, C(\hbox{pt}))$, which is identified with
relative $K^G$-homology. When $\overline M$ does not have finite $G$-homotopy
type, $K^G$-homology here means Steenrod $K^G$-homology, as explained in 
\cite{KKS}.} We will be restricting attention to the case where the
action of $G$ is proper and cocompact, in which case
$K^{-*}_G(C_0(M))$ may be viewed as a kind of orbifold $K$-homology
for the compact orbifold $G\backslash M$ (see \cite[Theorem 20.2.7]{Bl}.)

We will work throughout with real scalars and real $K$-theory, and
use a variant of the strategy found in \cite{Rosenberg(1999a)}
to prove Theorem \ref{the: Eul^G(X) and vector fields}.
\begin{proof}[Proof of Theorem \ref{the: Eul^G(X) and vector fields}]
Recall that since $\Xi$ is transverse to the zero-section, its
zero set $\Zero(\Xi)$ is discrete, and since $M$ is assumed
$G$-cocompact, $\Zero(\Xi)$ consists of only finitely many $G$-orbits.
Write $\Zero(\Xi)=\Zero(\Xi)^+\amalg \Zero(\Xi)^-$, according to the
signs of the indices $s(\Xi,x)$ of the zeros $x\in\Zero(\Xi)$. 
We fix a $G$-invariant
Riemannian metric on $M$ and use it to identify the form bundle of
$M$ with the Clifford algebra bundle $\Cl$ of the tangent bundle, 
with its standard grading in which vector fields are sections of
$\Cl^-$, and  
$D$ with the Dirac operator on $\Cl$.\footnote{Since sign conventions 
differ, we emphasize that for us, unit tangent vectors on $M$ have
square $-1$ 
in the Clifford algebra.} (This is legitimate by \cite[II, Theorem
5.12]{LM}.)  Let $\cH=\cH^+\oplus \cH^-$ be the 
$\bZ/2$-graded Hilbert space of $L^2$
sections of $\Cl$. Let $A$ be the operator
on $\cH$ defined by \emph{right} Clifford multiplication by $\Xi$ on
$\Cl^+$ (the even part of $\Cl$) and by right
Clifford multiplication by $-\Xi$ on $\Cl^-$ (the odd part). 
We use \emph{right} Clifford multiplication since it commutes
with the symbol of $D$. Observe that $A$ is self-adjoint, with square 
equal to multiplication by the non-negative function $\vert \Xi(x)\vert^2$. 
Furthermore, $A$ is odd with respect to the grading and 
commutes with multiplication by scalar-valued functions.

For $\lambda\ge 0$, let $D_\lambda=D+\lambda A$. 
As in \cite{Rosenberg(1999a)}, each $D_\lambda$ defines an unbounded
$G$-equivariant Kasparov module in 
the same Kasparov class as $D$. In the ``bounded picture''
of Kasparov theory, the corresponding operator is 
\begin{equation}
\label{eq:Blambda}
B_\lambda=D_\lambda\bigl(1+D_\lambda^2\bigr)^{-\frac12}
=\frac{1}{\lambda}D_\lambda\left(\frac{1}{\lambda^2}+
\frac{1}{\lambda^2}D_\lambda^2\right)^{-\frac12}.
\end{equation}
The axioms satisfied by this operator that insure that it defines
a Kasparov $K^G$-homology class (in the ``bounded picture'')
are the following:
\begin{description}
\label{desc:Kaspconds}
\item[(B1)] It is self-adjoint, of norm $\le 1$, and commutes with the
action of $G$.
\item[(B2)] It is odd with respect to the grading of $\Cl$.
\item[(B3)] For $f\in C_0(M)$, $fB_\lambda\sim B_\lambda f$ and 
$fB_\lambda^2\sim f$, where $\sim$ denotes equality modulo compact
operators. 
\end{description}
We should point out that
(B1) is somewhat stronger than it needs to be when $G$ is infinite.
In that case, we can replace invariance of $B_\lambda$ under $G$ by
``$G$-continuity,'' the requirement (see \cite{K1} and \cite[\S
20.2.1]{Bl}) that  
\begin{description}
\item[(B1$'$)] 
$f(g\cdot B_\lambda - B_\lambda) \sim 0$ for $f\in C_0(M)$, $g\in G$.
\label{desc:Kaspconds'}
\end{description}

In order to simplify the calculations that are coming next,
we may assume without loss of generality that we've chosen the
$G$-invariant Riemannian metric on $M$ so that for each $z\in
\Zero(\Xi)$, in some small open $G_z$-invariant neighborhood $U_z$ of
$z$, $M$ is $G_z$-equivariantly isometric to a ball, say of radius $1$, 
about the origin in Euclidean space
$\bR^n$ with an orthogonal $G_z$-action, with $z$ corresponding to the origin. 
This can be arranged since the exponential map induces a $G_z$-diffeomorphism of a small 
$G_z$-invariant neighborhood of $0 \in T_zM$ onto a $G_z$-invariant neighborhood of $z$
such that $0$ is mapped to $z$ and its differential at $0$ is the identity on $T_zM$ under the 
standard identification $T_0(T_zM) = T_zM$. 
Thus
the usual coordinates $x_1,\,x_2,\,\dots,\,x_n$ in Euclidean space
give local coordinates in $M$ for $|x|<1$, and
$\frac{\partial}{\partial x_1},\,\frac{\partial}{\partial
x_2},\,\dots,\,\frac{\partial}{\partial x_n}$
define a local orthonormal frame in $TM$ near $z$. We can arrange that
$(\rr^n)^{G_z}$ contains the points with $x_2 = \ldots = x_n = 0$ if
$(\rr^n)^{G_z}$ is different from $\{0\}$. In these exponential
local coordinates, the point $x_1=x_2=\cdots=x_n=0$ corresponds to $z$.
We may assume we
have chosen the vector field $\Xi$ so that in these local coordinates,
$\Xi$ is given by the radial vector field 
\begin{equation}
\label{eq:radialvecfield}
x_1\frac{\partial}{\partial x_1}+ x_2\frac{\partial}{\partial x_2}+ \cdots 
+x_n\frac{\partial}{\partial x_n}
\end{equation}
if $z\in \Zero(\Xi)^+$, or by the vector field 
\begin{equation}
\label{eq:negindexvecfield}
-x_1\frac{\partial}{\partial x_1}+ x_2\frac{\partial}{\partial x_2}+ \cdots 
+x_n\frac{\partial}{\partial x_n}
\end{equation}
if $z\in \Zero(\Xi)^-$. Thus $|\Xi(x)|=1$ on $\partial U_z$ for each $z$,
and we can assume (rescaling $\Xi$ if necessary) that $|\Xi|\ge 1$ on the
complement of $\bigcup_{z\in \Zero(\Xi)} U_z$. Recall that 
$D_\lambda=D+\lambda A$.
\begin{lem}
\label{lem:Dlambdalowerbound}
Fix a small number $\varepsilon>0$, and let $P_\lambda$ denote the
spectral projection of $D_\lambda^2$ corresponding to
$[0,\varepsilon]$. Then for $\lambda$ sufficiently large, 
$\range P_\lambda$ is $G$-isomorphic to $L^2(\Zero(\Xi))$ {\lp}a
Hilbert space with $\Zero(\Xi)$ as orthonormal basis, with the obvious
unitary action of $G$ coming from the action of $G$ on
$\Zero(\Xi)${\rp}, and there is a constant $C>0$ such that  
$(1-P_\lambda)D_\lambda^2\ge
C\lambda$. {\lp}In other words, $(\varepsilon,\,C\lambda)
\cap (\operatorname{spec}D_\lambda^2) = \emptyset$.{\rp}
Furthermore, the functions in
$\range P_\lambda$ become increasingly concentrated near $\Zero(\Xi)$ as
$\lambda\to \infty$. 
\end{lem}
\begin{proof}First observe that in Euclidean space $\bR^n$, if $\Xi$ is
defined by \eqref{eq:radialvecfield} or \eqref{eq:negindexvecfield}
and $A$ and $D_\lambda$ are defined from $\Xi$ as on $M$, then 
$S_\lambda=D_\lambda^2$ is
basically a Schr\"odinger operator for a harmonic oscillator, so one
can compute its spectral decomposition explicitly. (For example, if
$n=1$, then $S_\lambda=-\frac{d^2}{dx^2}+\lambda^2x^2\pm\lambda$,
the sign depending on whether $z\in \Zero(\Xi)^+$ or 
$z\in \Zero(\Xi)^-$ and whether
one considers the action on $\cH^+$ or $\cH^-$.)
When $z\in \Zero(\Xi)^+$, the kernel of $S_\lambda$
in $L^2$ sections of $\text{Cliff}(T\bR^n)$
is spanned by the Gaussian function
\[
(x_1,\,x_2,\,\dots,\,x_n)\mapsto e^{-\lambda \vert x\vert^2/2},
\]
and if $z\in \Zero(\Xi)^-$, the $L^2$ kernel is spanned by a similar
section of $\Cl^-\!\!$, $e^{-\lambda \vert x\vert^2/2}
\frac{\partial}{\partial x_1}$. Also, in both cases, $S_\lambda$ has
discrete spectrum lying on an arithmetic progression,
with one-dimen\-sion\-al kernel (in $L^2$)
and first non-zero eigenvalue given by $2n\lambda$.

Now let's go back to the operator on $M$.
Just as in \cite[Lemma 2]{Rosenberg(1999a)}, we have the estimate 
\begin{equation}
\label{eq:Dlambdasq}
-K\lambda \le D_\lambda^2-(D^2+\lambda^2 A^2)\le K\lambda,
\end{equation}
where $K>0$ is some constant (depending on the size of the covariant
derivatives of $\Xi$).\footnote{We are using the cocompactness of the
$G$-action to obtain a uniform estimate.} But $D_\lambda^2\ge 0$, and also,
from \eqref{eq:Dlambdasq},
\begin{equation}
\label{eq:Dlambdasqlowerbound}
\frac{1}{\lambda^2}D_\lambda^2\ge A^2 +
\frac{1}{\lambda^2}{D^2} - \frac{K}{\lambda},
\end{equation}
which implies that 
\begin{equation}
\label{eq:Dlambdasqlowerbound2}
\frac{1}{\lambda^2}D_\lambda^2\ge \text{multiplication by }
|\Xi(x)|^2-\frac{K}{\lambda}.
\end{equation}
So if $\xi_\lambda$ is a unit vector in $\range P_\lambda$, we have
\begin{equation}
\label{eq:Dlambdasqlowerbound3}
\frac{\varepsilon}{\lambda^2}\ge\left\langle
\frac{1}{\lambda^2}D_\lambda^2\xi_\lambda,\,\xi_\lambda 
\right\rangle\ge 
\int_M |\Xi(x)|^2\,|\xi_\lambda(x)|^2\,\dvol
-\frac{K}{\lambda}\bigl\Vert \xi_\lambda\bigr\Vert^2.
\end{equation}
Now $\Vert \xi_\lambda\Vert =1$, and if we fix $\eta>0$,
we only make the integral smaller by replacing $|\Xi(x)|^2$ by
$\eta$ on the set $E_\eta=\{x:|\Xi(x)|^2\ge \eta\}$ and by $0$
elsewhere. So
\[
\frac{\varepsilon}{\lambda^2} \ge -\frac{K}{\lambda}+ \eta
\int_{E_\eta} |\xi_\lambda(x)|^2\,\dvol
\]
or
\begin{equation}
\label{eq:kernelconc}
\left\Vert
\xi_\lambda\chi_{E_\eta}\right\Vert^2\le
\frac{K}{\eta\lambda} + \frac{\varepsilon}{\eta\lambda^2}.
\end{equation}
This being true for any $\eta$, we have verified
that as $\lambda\to\infty$, $\xi_\lambda$ becomes increasingly
concentrated near the zeros of $\Xi$, in the sense that the $L^2$ norm
of its restriction to the complement of any neighborhood of $\Zero(\Xi)$ goes
to $0$.

It remains to compute $\range P_\lambda$ (as a unitary representation
space of $G$) and to prove that
$D_\lambda^2$ has the desired spectral gap. Define a 
$C^2$ cut-off function
$\varphi(t)$, $0\le t<\infty$, so that $0\le \varphi(t)\le 1$,
$\varphi(t)=1$ for $0\le t\le \frac12$, $\varphi(t)=0$ for $t\ge 1$,
and $\varphi$ is decreasing on the interval
$\left[\frac12,\,1\right]$. In other words, $\varphi$ is supposed to
have a graph like this:
\begin{center}
\includegraphics[scale=.7]{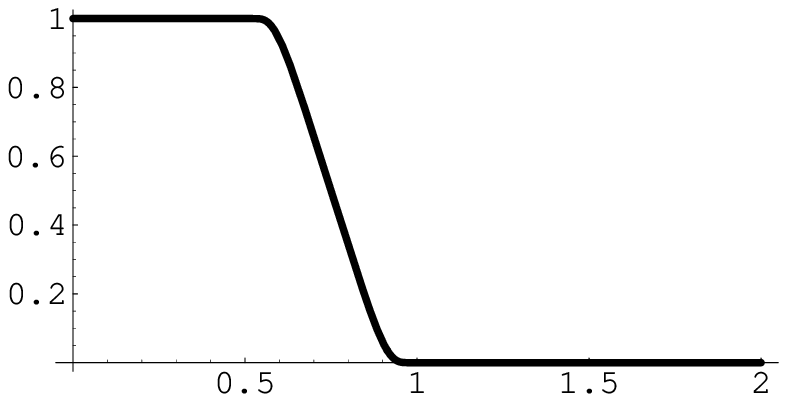}
\end{center}
We can arrange that $|\varphi'(t)|\le 3$ and that
$|\varphi''(t)|\le 20$.
For each element $z$ of $\Zero(\Xi)$, recall that we have a $G_z$-invariant neighborhood $U_z$
that can be identified with the unit ball in $\bR^n$ equipped with an orthogonal $G_z$-action.
So the function 
$\psi_{z,\lambda}(x)=\varphi(r)e^{-\lambda r^2/2}$, where $r=|x|$ is
the radial coordinate in 
$\bR^n$, makes sense as a function in $C^2(M)$, with support in
$U_z$. For simplicity suppose $z\in
\Zero(\Xi)^+$; the other case is exactly analogous except that we need a
$1$-form instead of a function. Then $D_\lambda^2$, acting
on radial functions, becomes
\[
-\Delta + \lambda^2|x|^2 -n\lambda = -\frac{\partial^2}{\partial r^2}
- (n-1)\frac{1}{r}\frac{\partial}{\partial r}+ \lambda^2r^2 - n\lambda.
\]
As we mentioned before, this operator on $\bR^n$ annihilates
$x\mapsto e^{-\lambda r^2/2}$, so we have
\begin{align}
\label{eq:Rayleighquotient}
\frac{\Vert D_\lambda \psi_{z,\lambda} \Vert^2}
{\Vert \psi_{z,\lambda} \Vert^2}
&=\frac{
\left \langle D_\lambda^2 \psi_{z,\lambda},\, 
 \psi_{z,\lambda}\right\rangle}
{\left \langle  \psi_{z,\lambda},\, \psi_{z,\lambda}\right\rangle} \notag\\
&=
\frac{
\int_0^1 \varphi(r)
\left( -r\varphi''(r) + (1 - n + 2 r^2 \lambda) \varphi'(r)\right)
e^{-\lambda r^2} r^{n-2}\,dr}
{\int_0^1 \varphi(r)^2e^{-\lambda r^2} r^{n-1}\,dr}
\notag\\
&\le
\frac{
\int_{1/2}^1 (20r+6\lambda)e^{-\lambda r^2} r^{n-2}\,dr}
{\int_0^{1/2} e^{-\lambda r^2} r^{n-1}\,dr}\,.
\end{align}
The expression \eqref{eq:Rayleighquotient} goes to $0$ faster than
$\lambda^{-k}$ for any $k\ge 1$, since the numerator dies rapidly and
the denominator behaves like a constant times
$\lambda^{-n/2}$ for large $\lambda$,
so $P_\lambda\psi_{z,\lambda}$ is non-zero and very close to
$\psi_{z,\lambda}$. Rescaling constructs a
unit vector in $\range P_\lambda$ concentrated near $z$, regardless of
the value of $\varepsilon$, provided $\lambda$ is sufficiently large
(depending on $\varepsilon$). And the action of $g\in G$ sends this
unit vector to the corresponding unit vector concentrated near $g\cdot z$.
In particular, $\range P_\lambda$ contains a Hilbert space
$G$-isomorphic to $L^2(\Zero(\Xi))$.

To complete the proof of the Lemma, it will suffice to show that if
$\xi$ is a unit vector in the domain of $D$
which is orthogonal to each $\psi_{z,\lambda}$, then $\Vert D_\lambda
\xi\Vert^2 \ge C\lambda$ for some constant $C>0$, provided $\lambda$ is sufficiently large 
Let $E=\bigcup_{z\in \Zero(\Xi)} V_z$, where $V_z$ corresponds to the ball
about the origin of radius $\frac12$ when we identify $U_z$ with the
ball about the origin in $\bR^n$ of radius $1$. Let $\chi_E$ be the
characteristic function of $E$. Then 
\[
1=\Vert \xi\Vert^2 = \Vert \chi_E\xi\Vert^2 + \Vert
(1-\chi_E)\xi\Vert^2. 
\]
Hence we must be in one of the following two cases:
\begin{enumerate}
\item $\Vert (1-\chi_E)\xi\Vert^2 \ge \frac12$.
\item $\Vert \chi_E\xi\Vert^2 \ge \frac12$.
\end{enumerate}
In case (1), we can argue just as in the inequalities
\eqref{eq:Dlambdasqlowerbound3} and \eqref{eq:kernelconc} with
$\eta=\frac{1}{4}$, since $E$ is precisely the set where
$|\Xi(x)|^2<\frac{1}{4}$. So we obtain 
\[
\frac{1}{\lambda^2} \Vert D_\lambda\xi\Vert^2
=\left\langle
\frac{1}{\lambda^2}D_\lambda^2\xi,\,\xi \right\rangle\ge 
-\frac{K}{\lambda} + \frac{1}{4}\left\Vert (1-\chi_E)\xi\right\Vert^2 
\ge \frac{1}{8}-\frac{K}{\lambda},
\]
which gives $\Vert D_\lambda \xi\Vert^2\succeq \text{const}\cdot\lambda^2$
once $\lambda$ is sufficiently large.  So now consider
case (2). Then for some $z$, we must have $\Vert \chi_{G\cdot V_z}
\xi\Vert^2 \ge \frac{1}{2|G\backslash\!
\Zero(\Xi)|}$. But by assumption, $\xi \perp \psi_{g\cdot z,\lambda}$
(for this same $z$ and all $g\in G$). 
Assume for simplicity that $\xi\in\cH^+$ and
$z\in \Zero(\Xi)^+$.  If $\xi\in\cH^+$ and
$z\in \Zero(\Xi)^-$, there is no essential difference,
and if $\xi\in\cH^-$, the calculations are similar, but we need $1$-forms
in place of functions.  Anyway, if we let $\xi_g$ denote
$\xi|_{U_g\cdot z}$ transported to $\bR^n$, we have 
\[
\begin{aligned}
0=&\int_{\bR^n} \varphi(|x|)\xi_g(x)e^{-\lambda|x|^2/2}\,dx
\qquad (\forall g),\\
1\ge
&\sum_g \int_{|x|\le\frac12}\varphi(|x|)^2 \bigl|\xi_g(x)\bigr|^2 
\, dx\ge \frac{1}{2|G\backslash\!\Zero(\Xi)|}.
\end{aligned}
\]
Now we use the fact that the Schr{\"o}dinger operator $S_\lambda$ on
$\bR^n$ has one-dimen\-sion\-al kernel in $L^2$
spanned by $x\mapsto e^{-\lambda|x|^2/2}$ (if $z\in \Zero(\Xi)^+$), and spectrum 
bounded below by $2n\lambda$ on the orthogonal complement of this
kernel. (If $z\in \Zero(\Xi)^-$, the entire spectrum of $S_\lambda$ on
$\cH^+$ is bounded below by $2n\lambda$.) So compute as follows:
\begin{align}
\left\Vert D_\lambda \bigl(\varphi(|x|)\xi_g\bigr)\right\Vert^2 &= 
\left\langle  D_\lambda^2\bigl(\varphi(|x|)\xi_g\bigr),\,
\varphi(|x|)\xi_g \right\rangle \notag\\
&\ge 2n\lambda \left\langle  \varphi(|x|)\xi_g,\,
\varphi(|x|)\xi_g \right\rangle .
\label{eq:Schrlowerbound}
\end{align}
Let $\omega$ be the function on $M$ which is $0$ on the complement of 
$\bigcup_{g} U_{g\cdot z}$ and given by $\varphi(|x|)$ on $U_{g\cdot
z}$ (when we use the local coordinate system there centered at $g\cdot z$).
Then:
\begin{align}
\left\Vert D_\lambda \xi\right\Vert^2 &= \left\Vert D_\lambda
\bigl(\omega\xi\bigr)\right\Vert^2 +
\left\Vert D_\lambda
\bigl(\bigl(1-\omega\bigr)\xi\bigr)\right\Vert^2 \notag\\
\label{eq:innprodterm}
&\qquad
+2\left\langle  D_\lambda^2\bigl(\bigl(1-\omega\bigr)\xi\bigr),\,
\omega\xi \right\rangle.
\end{align}
Since $D_\lambda$ is local and $\omega$ is supported on the $U_{g\cdot
z}$, $g\in G$, the first term on the right is simply
\begin{equation}
\label{eq:sumoverg}
\left\Vert D_\lambda
\bigl(\omega\xi\bigr)\right\Vert^2 =
\sum_g \left\Vert D_\lambda \bigl(\varphi(|x|)\xi_g\bigr)\right\Vert^2
\ge \frac{2n\lambda}{2|G\backslash\!\Zero(\Xi)|}
\end{equation}
by \eqref{eq:Schrlowerbound}.
In the inner product term in \eqref{eq:innprodterm}, since
$\omega\xi$ is a sum of pieces with disjoint supports $U_{g\cdot z}$,
we can split this as a sum over terms we can transfer 
to $\bR^n$, getting
\[
\begin{aligned}
2&\sum_g
\left\langle  S_\lambda\bigl(\bigl(1-\varphi(|x|)\bigr)\xi_g\bigr),\,
\varphi(|x|)\xi_g \right\rangle
 = 2\sum_g\left\langle D^2\bigl(\bigl(1-\varphi(|x|)\bigr)\xi_g\bigr),\,
\varphi(|x|)\xi_g \right\rangle\\
&\qquad +2\sum_g\int_{\frac12 \le |x| \le 1} (\lambda^2|x|^2+T\lambda)\,
\varphi(|x|)
\bigl(1-\varphi(|x|)\bigr)|\xi_g(x)|^2\,dx,\\
\end{aligned}
\]
where $-K \le T\le K$. Since $\lambda^2|x|^2+T\lambda>0$ on
$\frac12 \le |x| \le 1$ for large enough $\lambda$, the integral here
is nonnegative, and the only possible negative contributions to 
$\left\Vert D_\lambda \xi\right\Vert^2$ are the terms
$ 2\left\langle D^2\bigl(\bigl(1-\varphi(|x|)\bigr)\xi_g\bigr),\,
\varphi(|x|)\xi_g \right\rangle$,
which do not grow with $\lambda$. So from
\eqref{eq:innprodterm}, \eqref{eq:sumoverg}, and
\eqref{eq:Schrlowerbound},  $\left\Vert D_\lambda \xi\right\Vert^2 \ge 
\text{const}\cdot\lambda$ for large enough $\lambda$, which
completes the proof. 
\end{proof}

\textbf{Proof of Theorem \ref{the: Eul^G(X) and vector fields},
continued}\qua
We begin by defining a continuous field $\cE$ of $\bZ/2$-graded
Hilbert spaces over  
the closed interval $[0,\,+\infty]$.  Over the open interval 
$[0,\,+\infty)$, the field is just the trivial one, with fiber
$\cE_\lambda=\cH$, the $L^2$ sections of $\Cl$. But the fiber
$\cE_\infty$ over $+\infty$ will 
be the direct sum of $\cH\oplus V$, where $V=L^2(\Zero(\Xi))$ is a
Hilbert space with orthonormal basis $v_z$, $z\in \Zero(\Xi)$. We put a
$\bZ/2$-grading on $V$ by letting $V^+=L^2(\Zero(\Xi)^+)$,
$V^-=L^2(\Zero(\Xi)^-)$. To define the continuous field 
structure, it is enough by \cite[Proposition 10.2.3]{Dix} to define a
suitable total set of continuous sections near the exceptional point
$\lambda=\infty$. We will declare ordinary continuous functions
$[0,\infty]\to \cH$ to be continuous, but will also allow
additional continuous fields that become increasingly concentrated
near the points of $\Zero(\Xi)$. Namely, suppose $z\in \Zero(\Xi)$. By Lemma
\ref{lem:Dlambdalowerbound}, for $\lambda$
large, $D_\lambda$ has an element $\psi_{z,\lambda}$ 
in its ``approximate kernel'' increasingly
supported close to $z$, and we have a formula for it. So we declare
$(\xi(\lambda))_{\lambda<\infty}$ to define a continuous field converging
to $cv_z$ at $\lambda=\infty$ if for any neighborhood $U$ of $z$,
\[
\int_{M\smallsetminus U} |\xi(\lambda)(m)|^2\,\dvol(m)\to 0
\qquad\text{as }\lambda\to\infty,
\]
and if (assuming $z\in \Zero(\Xi)^+$) $\xi(\lambda)\in\cH^+$ and 
\begin{equation}
\label{eq:convtovz}
\left\Vert \xi(\lambda) - c \left(\frac{\lambda}{\pi}\right)^{\frac
n4} \psi_{z,\lambda}
\right\Vert \to 0\qquad\text{as }\lambda\to\infty.
\end{equation}
The constant reflects the fact that the $L^2$-norm of
$e^{-\lambda \vert x\vert^2/2}$ is $\left(\frac{\pi}{\lambda}\right)^{\frac
n4}$.  If $z\in \Zero(\Xi)^-$, we use the same definition, but require
$\xi(\lambda)\in\cH^-$.

This concludes the definition of the continuous field of Hilbert
spaces $\cE$, which we can think of as a Hilbert $C^*$-module over
$C(I)$, $I$ the interval $[0,\,+\infty]$. We will use this to define a
Kasparov $(C_0(M),\,C(I))$-bimodule, or in other words, a homotopy
of Kasparov $(C_0(M),\,\bR)$-modules. The action of $C_0(M)$ on $\cE$ is
the obvious one: $C_0(M)$ acts on $\cH$ the usual way, and it acts on
$V$ (the other summand of $\cE_\infty$) by evaluation of functions at
the points of $\Zero(\Xi)$: 
\[
f\cdot v_z = f(z)v_z,\qquad z\in \Zero(\Xi),\,\,f\in C_0(M).
\]
We define a field $T$ of operators  on $\cE$ as follows.
For $\lambda<\infty$, $T_\lambda\in\cL(\cE_\lambda)=\cL(\cH)$ is simply
$B_\lambda$ as defined in \eqref{eq:Blambda}, where recall that
$D_\lambda=D+\lambda A$. For $\lambda=\infty$, $\cE_\infty=\cH\oplus
V$ and $T_\infty$ is $0$ on $V$ and is given on $\cH$ by 
the operator $B_\infty=A/\vert A\vert$ which is right
Clifford multiplication by $\frac{\Xi(x)}{\vert \Xi(x)\vert}$ 
(an $L^\infty$, but possibly discontinuous, vector field) on $\cH^+$
and by $\frac{-\Xi(x)}{\vert \Xi(x)\vert}$ on $\cH^-$. Note that
$T_\infty^2$ is $0$ on $V$ and the identity on $\cH$. While
$1_V$ is not compact on $V$ if $\Zero(\Xi)$ is infinite, this is not a
problem since for $f\in C_c(M)$, the action of $f$ on $V$ has finite
rank (since $f$ annihilates $v_z$ for $z\not\in \supp f$).

Now we check the axioms for $(\cE,\,T)$ to define a homotopy of
Kasparov modules from $[D]$ to the class of 
\[
\bigl( C_0(M),\, \cE_\infty,\,T_\infty \bigr)= 
\bigl( C_0(M),\, \cH,\,B_\infty \bigr) \oplus 
\bigl( C_0(M),\, V,\,0 \bigr) . 
\]
But $\bigl( C_0(M),\,
\cH,\,B_\infty \bigr)$ is a \emph{degenerate} Kasparov module,
since $B_\infty$ commutes with multiplication by functions and
has square $1$. So the class of $\bigl( C_0(M),\,
\cE_\infty$, $T_\infty \bigr)$ is just the class of $\bigl( C_0(M),\,
V,\,0 \bigr)$, which (essentially by definition)
is the image under the inclusion
$\Zero(\Xi)\hookrightarrow M$ of the sum (over $G\backslash\!
\Zero(\Xi)$) of $+1$ times the canonical
class $KO^G_0(z)\circ (\alpha_z)_*([\bR])$
for $G\cdot z\subseteq \Zero(\Xi)^+$ and of
$-1$ times this class if $G\cdot z\subseteq \Zero(\Xi)^-$. 
This will establish Theorem \ref{the: Eul^G(X) and vector fields},
assuming we can verify that we have a homotopy of Kasparov modules.

The first thing to check is that the action of $C_0(M)$ on $\cE$ is
continuous, i.e., given by a $*$-homomorphism $C_0(M)\to\cL(\cE)$. The
only issue is continuity at $\lambda=\infty$. In other words, since
the action on $\cH$ is constant, we just need to know that if $\xi$ is
a continuous field converging as $\lambda\to \infty$ to a vector $v$ in
$V$, then for $f\in C_0(M)$, $f\cdot \xi(\lambda)\to f\cdot v$. But it's
enough to consider the special kinds of continuous fields discussed
above, since they generate the structure, and if $\xi(\lambda)\to c
v_z$, then $\xi(\lambda)$ becomes increasingly concentrated at $z$
(in the sense of $L^2$ norm), and hence $f\cdot  \xi(\lambda) \to c
f(z)v_z$, as required. 

Next, we need to check that $T\in\cL(\cE)$. Again, the only issue is
(strong operator) continuity at $\lambda=\infty$. Because of the way
continuous fields are defined at $\lambda=\infty$, there are basically
two cases to check.  First, if $\xi\in\cH$, we need to check that
$B_\lambda \xi \to B_\infty \xi$ as $\lambda\to \infty$. Since the
$B_\lambda $'s all have norm $\le 1$, we also only need to check this
on a dense set of $\xi$'s. First, fix $\varepsilon >0$ small
and suppose $\xi$ is smooth and supported on the open set
where $|\Xi(x)|^2 > \varepsilon$. Then for $\lambda$ large,
Lemma \ref{lem:Dlambdalowerbound} implies that there is a constant
$C>0$ (depending on $\varepsilon$ but not on $\lambda$) such that
$\langle D_\lambda^2 \xi,\,\xi\rangle > C\lambda\Vert\xi\Vert^2$.
In fact, if $P_\lambda$ is the spectral projection of $D_\lambda^2$
for the interval $[0,\,C\lambda]$, 
Lemma \ref{lem:Dlambdalowerbound} implies that 
$\Vert P_\lambda\xi\Vert \le \varepsilon \Vert \xi\Vert $
for $\lambda$ sufficiently large.
(This is because the condition on the support of $\xi$ forces
$\xi$ to be almost orthogonal to the spectral subspace where
$ D_\lambda^2 \le \varepsilon$.) Now let $E^+_\lambda$ 
and $E^-_\lambda$ be the
spectral projections for $D_\lambda$ corresponding to the intervals
$(0,\,\infty)$ and $(-\infty,\,0)$, and let $F^+$ and $F^-$ be the
spectral projections for $A$ corresponding to the same intervals.
Since the vector field $\Xi$ vanishes only on a discrete set, the operator
$A$ has no kernel, and hence $F^+ + F^- = 1$. Now we appeal to
two results in Chapter VIII of \cite{Kato}:
Corollary 1.6 in \S1, and Theorem 1.15 in \S2.  The former shows
that the operators $A+\frac{1}{\lambda}D$, all defined on $\dom D$,
``converge strongly in the generalized sense'' to $A$. Since
the positive and negative spectral subspaces for
$A+\frac{1}{\lambda}D$ are the same as for $D_\lambda$ (since
the operators only differ by a homothety), \cite[Chapter VIII,
\S2, Theorem 1.15]{Kato} then
shows that $E^+_\lambda \to F^+$ and $E^-_\lambda \to F^+$
in the strong operator topology.  Note that the fact that $A$ has no
kernel is needed in these results.

Now since $\Vert P_\lambda\xi\Vert \le \varepsilon \Vert \xi\Vert $
for $\lambda$ sufficiently large, we also have
\[
B_\infty\xi = F^+\xi - F^-\xi,
\quad \text{and} \quad 
\Vert B_\lambda\xi - (E_\lambda^+\xi - E_\lambda^-\xi) \Vert 
\le 2\varepsilon
\]
for $\lambda$ sufficiently large. Hence
\[
\Vert  B_\lambda\xi - B_\infty\xi \Vert \le 2\varepsilon
+ \Vert  (E_\lambda^+\xi - F^+\xi) - (E_\lambda^-\xi - F^-\xi)
\Vert  \to 2\varepsilon.
\]
Now let $\varepsilon\to 0$. Since, with $\varepsilon$ tending to zero,
$\xi$'s satisfying our support condition are dense,
we have the required strong convergence.

There is one other case to check, that where $\xi(\lambda)\to c
v_z$ in the sense of the continuous field structure of $\cE$.
In this case, we need to show that $B_\lambda \xi(\lambda)\to 0$.
This case is much easier: $\xi(\lambda)\to c v_z$ means
\[
\left\Vert \xi(\lambda) - c \left(\frac{\lambda}{\pi}\right)^{\frac
n4} \psi_{z,\lambda}
\right\Vert \to 0\qquad \text{by \eqref{eq:convtovz}},
\]
while $\Vert B_\lambda\Vert\le 1$ and
\[
\left\Vert D_\lambda \left(\left(\frac{\lambda}{\pi}\right)^{\frac
n4} \psi_{z,\lambda}\right) \right\Vert \to 0
\qquad \text{by \eqref{eq:Rayleighquotient}},
\]
so $B_\lambda \xi(\lambda)\to 0$ in norm.

Thus $T\in\cL(\cE)$. Obviously, $T$ satisfies (B1) and (B2) of
page \pageref{desc:Kaspconds}, so we need to check the analogues of
(B3), which are that $f(1-T^2)$ and $[T,f]$ lie in $\cK(\cE)$ for $f\in
C_0(M)$. First consider $1-T^2$. $1-T_\lambda^2$ is locally compact
(i.e., compact after multiplying by $f\in C_c(M)$)
for each $\lambda$, since 
\[
1-T_\lambda^2=1-B_\lambda^2=(1+D_\lambda^2)^{-1}=
\bigl(1+ (D+\lambda A)^2 \bigr)^{-1}
\]
is locally compact for $\lambda<\infty$, and $1-T_\infty^2$ is just projection
onto $V$, where functions $f$ of compact support act by finite-rank
operators. So we just need to check that
$1-T_\lambda^2$ is a norm-continuous field of operators on $\cE$.
Continuity for $\lambda<\infty$ is routine, and implicit
in \cite[Remarques 2.5]{BJ}. To check continuity at $\lambda=\infty$,
we use Lemma \ref{lem:Dlambdalowerbound}, which shows that
$(1+D_\lambda^2)^{-1}=P_\lambda+O\bigl(\frac{1}{\lambda}\bigr)$,
and also that $P_\lambda$ is increasingly concentrated near $\Zero(\Xi)$. 
So near $\lambda=\infty$, we can write the field of operators
$(1+D_\lambda^2)^{-1}$  as a sum of rank-one projections onto vector fields
converging to the various $v_z$'s (in the sense of our continuous
field structure) and another locally compact operator converging in norm to
$0$.  

This leaves just one more thing to check, that for $f \in C_0(M)$, $[f,\,
T_\lambda]$ lies in $\cK(\cE)$. We already know that
$[f,\,B_\lambda]\in \cK(\cH)$ for fixed $\lambda$
and is norm-continuous in $\lambda$ for $\lambda<\infty$, so since
$T_\infty$ commutes with multiplication operators, it suffices to show
that $[f,\,B_\lambda]$ converges to $0$ in norm as $\lambda\to 0$. We
follow the method of proof in \cite[p.\ 3473]{Rosenberg(1999a)},
pointing out the 
changes needed because of the zeros of the vector field $\Xi$.

We can take $f\in C_c^\infty(M)$ with critical points at all of the points
of the set $\Zero(\Xi)$, since such functions are dense in $C_0(M)$.
Then estimate as follows:
\begin{align}
\label{eq:commutator}
[f,\,B_\lambda] &= \left[ f,\,D_\lambda(1+ D_\lambda^2)^{-1/2} \right]
\notag \\
& =[f,\,D_\lambda](1+ D_\lambda^2)^{-1/2} 
 + D_\lambda \left[f,\,(1+ D_\lambda^2)^{-1/2}\right].
\end{align}
We have $[f,\,D_\lambda]=[f,\,D]$, which is a $0$'th order operator
determined by the derivatives of $f$, of compact support since
$f$ has compact support, and we've seen that
$(1+ D_\lambda^2)^{-1/2}$ converges as $\lambda\to\infty$ (in the norm
of our continuous field) to 
projection onto the  space $V=L^2(\Zero(\Xi))$. 
Since the derivatives of $f$ vanish on $\Zero(\Xi)$, the product
$[f,\,D_\lambda](1+ D_\lambda^2)^{-1/2}$, which is
the first term in \eqref{eq:commutator},
goes to $0$ in norm. As for the second
term, we have (following \cite[p.\ 199]{Bl})
\begin{equation}
\label{eq:intformula}
D_\lambda \left[f,\,(1+ D_\lambda^2)^{-1/2}\right]
= \frac{1}{\pi}\int_0^\infty \mu^{-\frac12}D_\lambda
\left[f,\,(1+ D_\lambda^2+\mu)^{-1}\right] \,d\mu,
\end{equation}
and
\[
D_\lambda \left[f,\,(1+ D_\lambda^2+\mu)^{-1}\right] 
= D_\lambda (1+ D_\lambda^2+\mu)^{-1} \left[1+ D_\lambda^2+\mu,\,f\right]
(1+ D_\lambda^2+\mu)^{-1}.
\]
Now use the fact that 
\[
\left[1+ D_\lambda^2+\mu,\,f\right]=\left[D_\lambda^2,\,f\right]
=D_\lambda [D_\lambda,\,f] + [D_\lambda,\,f] D_\lambda = D_\lambda [D,\,f]
+[D,\,f] D_\lambda.
\]
We obtain that
\begin{multline}
\label{eq:commutatorintegrand}
D_\lambda \left[f,\,(1+ D_\lambda^2+\mu)^{-1}\right] \\
=\frac{D_\lambda^2}{1+D_\lambda^2+\mu} [D,\,f] \frac{1}{1+D_\lambda^2+\mu}
+ \frac{D_\lambda}{1+D_\lambda^2
+\mu} [D,\,f] \frac{D_\lambda}{1+D_\lambda^2+\mu}.
\end{multline}
Again a slight modification of the argument in \cite[p.\
3473]{Rosenberg(1999a)} 
is needed, since $D_\lambda$ has an
``approximate kernel'' concentrated near the
points of $\Zero(\Xi)$. So we estimate the norm of the right side of
\eqref{eq:commutatorintegrand} as follows:
\begin{align}
\left\Vert D_\lambda \left[f,\,(1+ D_\lambda^2+\mu)^{-1}\right]
\right\Vert & \le \left\Vert 
\frac{D_\lambda^2}{1+D_\lambda^2+\mu} [D,\,f]
\frac{1}{1+D_\lambda^2+\mu}
\right\Vert \label{eq:commutator1} \\
& + \left\Vert 
\frac{D_\lambda}{1+D_\lambda^2
+\mu} [D,\,f] \frac{D_\lambda}{1+D_\lambda^2+\mu}
\right\Vert. \label{eq:commutator2} 
\end{align}
The first term, \eqref{eq:commutator1}, is bounded by the second,
\eqref{eq:commutator2}, plus an additional commutator term:
\begin{equation}
\label{eq:commutator3}
\left\Vert 
\frac{D_\lambda}{1+D_\lambda^2
+\mu} [D, [D, f]] \frac{1}{1+D_\lambda^2+\mu}
\right\Vert.
\end{equation}
Now the contribution of the term \eqref{eq:commutator2} is
estimated by observing that the function
\[
\frac{x}{1+x^2+\mu},\qquad -\infty < x < \infty
\]
has maximum value $\frac{1}{2\sqrt{1+\mu}}$ at $x=\sqrt{1+\mu}$,
is increasing for $0<x<\sqrt{1+\mu}$,
and is decreasing to $0$ for $x>\sqrt{1+\mu}$.
Fix $\varepsilon>0$ small.
Since, by Lemma \ref{lem:Dlambdalowerbound}, $|D_\lambda|$ has spectrum
contained in $[0,\,\sqrt\varepsilon]
\cup [\sqrt{C\lambda},\,\infty)$, we find that
\begin{equation}
\label{eq:maxest}
\left\Vert \frac{D_\lambda}{1+D_\lambda^2 +\mu} \right\Vert \le
\begin{cases}
\frac{1}{2\sqrt{1+\mu}}, & \\
\qquad\sqrt{C\lambda} \le \sqrt{1+\mu},
\text{ or }\mu \ge C\lambda - 1,\\ & \\
\max\left(\frac{\sqrt{\varepsilon}}{1+\varepsilon+\mu},
\frac{\sqrt{C\lambda}}{1 + \mu + C\lambda}\right), &\\
\qquad\sqrt{C\lambda} \ge
\sqrt{1+\mu},
\text{ or }0\le \mu \le C\lambda - 1.&
\end{cases}
\end{equation}
Thus the contribution of the term \eqref{eq:commutator2}
to the integral in \eqref{eq:intformula} is
bounded by
\begin{align}
\label{eq:intest}
&\frac{\Vert [D,f]\Vert}{\pi}\int_0^\infty 
\left\Vert \frac{D_\lambda}{1+D_\lambda^2 +\mu} \right\Vert^2
\frac{1}{\sqrt{\mu}} \, d\mu \notag\\
&\le\frac{\Vert [D,f]\Vert}{\pi}
\left(\int_{C\lambda - 1}^\infty \frac{1}{\sqrt\mu}
\frac{1}{4(1+\mu)}\,d\mu \right.\\
&\qquad + \left.
\int_0^{C\lambda - 1} \frac{1}{\sqrt\mu}\,\max\left(
\frac{C\lambda}{(1 + \mu + C\lambda)^2},\,
\frac{\varepsilon}{(1+\varepsilon+\mu)^2}\right)\,d\mu \right) \notag\\
&\le \frac{\Vert [D,f]\Vert}{\pi}
\left(\frac{1}{4}\int_{C\lambda-1}^\infty \mu^{-\frac32} \,d\mu
+ \int_0^{C\lambda} \frac{1}{\sqrt\mu}\,\frac{C\lambda}{(C\lambda)^2}  
\,d\mu + \int_0^\infty \frac{\varepsilon}{\sqrt{\mu}
(1+\mu)^2} \,d\mu \right)\notag\\
&= \frac{\Vert [D,f]\Vert}{\pi}
\left( \frac{1}{2\sqrt{C\lambda -1}} +
\frac{2}{\sqrt{C\lambda}} + \frac{\pi\varepsilon}{2}\right)
\to \frac{\Vert [D,f]\Vert}{2} \varepsilon.
\end{align}
We can make this as small as we like by taking $\varepsilon$ small enough.
Similarly, the contribution of term \eqref{eq:commutator3}
to the integral in \eqref{eq:intformula} is
bounded by
\begin{align}
\label{eq:intest2}
&\frac{\Vert [D,[D,f]]\Vert}{\pi}\int_0^\infty 
\left\Vert \frac{D_\lambda}{1+D_\lambda^2 +\mu} \right\Vert
\frac{1}{1+\mu}\,\frac{1}{\sqrt{\mu}}\,d\mu \notag\\
&\le \frac{\Vert [D,[D,f]]\Vert}{\pi}
\Bigl(\int_{C\lambda - 1}^\infty \frac{1}{2\sqrt{1+\mu}}\,\frac{1}{1+\mu}\,
\frac{1}{\sqrt{\mu}}\,d\mu \notag \\
& \qquad + 
\int_0^{C\lambda - 1} 
\frac{\sqrt{C\lambda}}{(1 + \mu + C\lambda)}\,\frac{1}{1+\mu}\,
\frac{1}{\sqrt{\mu}} \,d\mu  + \int_0^\infty \frac{\varepsilon}{\sqrt{\mu}
(1+\mu)^2} \,d\mu \Bigr) \notag\\
&\le \frac{\Vert [D,[D,f]]\Vert}{\pi}
\left( \int_{C\lambda - 1}^\infty \frac{1}{2\mu^2}\,d\mu
+ \int_0^\infty \frac{1}{\sqrt{C\lambda}}\,\frac{1}{{\sqrt\mu}(1+\mu)}
\,d\mu \right. \notag \\
& \qquad \qquad \qquad \qquad \left. 
+ \int_0^\infty \frac{\varepsilon}{\sqrt{\mu} (1+\mu)^2} \,d\mu 
\right) \notag \\
&\le \frac{\Vert [D,[D,f]]\Vert}{\pi} \left( 
\frac{1}{2(C\lambda-1)} + \frac{\pi}{\sqrt{C\lambda}}
+ \frac{\pi\varepsilon}{2}\right) \to \frac{\Vert [D,[D,f]]\Vert}{2}
\varepsilon,
\end{align}
which again can be taken as small as we like.
This completes the proof.
\end{proof}


\tit{Review of notions of equivariant Euler characteristic}
{Review of notions of equivariant Euler characteristic}

Next we briefly review the universal equivariant Euler characteristic,
as well as some other notions of equivariant Euler characteristic,
so we can see exactly how they are related to the $KO^G$-Euler
class $\Eul^G(M)$. We will use the following notation in the sequel.

\begin{notation} \label{not: X^H(x) and so on}
Let $G$ be a discrete group and $H \subseteq G$ be a subgroup.
Let $N\!H = \{g \in G \mid gHg^{-1} = H\}$ be its \emph{normalizer}
and let $W\!H := N\!H/H$ be its \emph{Weyl group}.

Denote by $\consub(G)$ the set of conjugacy classes $(H)$ of subgroups
$H \subseteq G$. 

Let $X$ be a $G$-$CW$-complex. Put 
\begin{eqnarray*}
X^H & := & \{x \in X \mid H \subseteq G_x\};
\\
X^{>H} & := & \{x \in X \mid H \subsetneq G_x\},
\end{eqnarray*} 
where $G_x$ is the isotropy group of $x$ under the $G$-action.

Let $x\co G/H \to X$ be a $G$-map. Let $X^H(x)$ be the component of
$X^H$ containing $x(1H)$. Put 
$$X^{>H}(x) = X^H(x) \cap X^{>H}.$$
Let $W\!H_x$ be the isotropy group of $X^H(x) \in \pi_0(X^H)$ under
the $W\!H$-action. 
\end{notation}

Next we define the group $U^G(X)$, in which the universal equivariant Euler characteristic
takes its values.
Let $\Pi_0(G,X)$ be the \emph{component category} of the $G$-space $X$
in the sense of  tom Dieck \cite[I.10.3]{Dieck(1987)}.
Objects are $G$-maps $x\co G/H \to X$. A morphism $\sigma$ from
$x\co G/H \to X$ to 
$y\co G/K \to X$ is a $G$-map $\sigma\co G/H \to G/K$ 
such that $y \circ \sigma$ and $x$ are $G$-homotopic. A $G$-map
$f\co  X \to Y$ induces a functor $\Pi_0(G,f)\co  \Pi_0(G,X) \to
\Pi_0(G,Y)$ 
by composition with $f$. Denote by $\Is \Pi_0(G,X)$ the set of
isomorphism classes $[x]$ of objects $x\co G/H \to X$ in
$\Pi_0(G,X)$. Define  
\begin{eqnarray}
U^G(X) & := & \zz[\Is \Pi_0(G,X)], 
\label{definition of U^G(X)}
\end{eqnarray}
where for a set $S$ we denote by $\zz [S]$ 
the free abelian group with basis $S$.
Thus we obtain a covariant functor from the category of $G$-spaces to the 
category of abelian groups. Obviously $U^G(f) = U^G(g)$ if $f,g \co
X \to Y$ are $G$-homotopic. 

There is a natural bijection
\begin{eqnarray}
\Is \Pi_0(G,X)  & \xrightarrow{\cong} & \coprod_{(H) \in \consub(G)} 
W\!H\backslash \pi_0(X^H),
\label{identifying Is Pi_0(G,X)}
\end{eqnarray}
which sends $x\co G/H \to X$ to the orbit under the $W\!H$-action on $\pi_0(X^H)$
of the component $X^H(x)$ of $X^H$ which contains the point $x(1H)$.
It induces a natural isomorphism
\begin{eqnarray}
U^G(X)  & \xrightarrow{\cong} & \bigoplus_{(H) \in \consub(G)} ~
\bigoplus_{W\!H\backslash \pi_0(X^H)} ~ \zz. 
\label{identifying U^G(X)}
\end{eqnarray}

\begin{defn} 
\label{def: universal equivariant Euler characteristic}
Let $X$ be a finite $G$-$CW$-complex $X$. We define the
\emph{universal equivariant Euler characteristic} of $X$
\begin{eqnarray*}
\chi^G(X) & \in & U^G(X)
\end{eqnarray*}
by assigning to $[x\co G/H \to X] \in \Is \Pi_0(G,X)$ the
(ordinary) Euler characteristic of the pair of finite $CW$-complexes
$(W\!H_x\backslash X^H(x),W\!H_x\backslash X^{>H}(x))$.

If the action of $G$ on $X$ is proper (so that the isotropy group of
any open cell in $X$ is finite), we define the \emph{orbifold Euler
characteristic} of $X$ by:
\begin{eqnarray*}
\chi^{\qq G}(X) & := & \sum_{p \ge 0} ~ \sum_{G\cdot e \in G\backslash
I_p(X)} |G_e|^{-1} 
\hspace{5mm} \in \qq,
\end{eqnarray*}
where $I_p(X)$ is the set of open cells of $X$ (after forgetting the
group action). 
\end{defn}

The orbifold Euler characteristic $\chi^{\qq G}(X)$ can be identified with
the more general notion of the $L^2$-Euler characteristic
$\chi^{(2)}(X;\caln(G))$, where $\caln(G)$ is the group von Neumann
algebra of $G$. One can compute $\chi^{(2)}(X;\caln(G))$ in
terms of $L^2$-homology 
$$\chi^{(2)}(X;\caln(G)) ~ = ~ 
\sum_{p \ge 0} (-1)^p \cdot
\dim_{\caln(G)}\left(H_p^{(2)}(X;\caln(G)\right),$$ 
where $\dim_{\caln(G)}$ denotes the von Neumann dimension (see for instance
\cite[Section 6.6]{Lueck(2002)}).

Next we define for a proper $G$-$CW$-complex $X$ the \emph{character map}
\begin{eqnarray}
\ch^G(X)\co U^G(X) & \to & \bigoplus_{\Is \Pi_0(G,X)}\qq ~ = ~ 
\bigoplus_{(H)\in \consub(G)} ~ \bigoplus_{W\!H\backslash \pi_0(X^H)} \qq.
\label{definition of character map chi^G(X)}
\end{eqnarray}
We have to define for an isomorphism class $[x]$ of objects $x\co G/H
\to X$ in   
$\Pi_0(G,X)$ the component $\ch^G(X)([x])_{[y]}$ of $\ch^G(X)([x])$
which belongs to an isomorphism class $[y]$ of objects $y\co G/K \to X$ in 
$\Pi_0(G,X)$, and check that $\chi^G(X)([x])_{[y]}$ is different from zero for
at most finitely many $[y]$. Denote by $\mor(y,x)$ the set of
morphisms from $y$ to $x$ in 
$\Pi_0(G,X)$. We have the left operation 
$$\aut(y,y) \times \mor(y,x) \to \mor(y,x), 
\hspace{5mm}
(\sigma, \tau) \mapsto \tau \circ \sigma^{-1}.$$
There is an isomorphism of groups
$$W\!K_y \xrightarrow{\cong} \aut(y,y)$$
which sends $gK \in W\!K_y$ to  the automorphism of  $y$ given by the
$G$-map $$R_{g^{-1}}\co G/K \to G/K, 
\hspace{5mm}
g'K \mapsto g'g^{-1}K.$$
Thus $\mor(y,x)$ becomes a left $W\!K_y$-set. 

The $W\!K_y$-set $\mor(y,x)$ can be rewritten as
$$\mor(y,x) ~ = ~ \{g \in G/H^K \mid g\cdot x(1H) \in X^K(y)\},$$
where the left operation of $W\!K_y$ on 
$\{g \in G/H^K \mid g\cdot x(1H) \in Y^K(y)\}$ comes from the canonical
left action of $G$ on $G/H$. Since $H$ is finite and hence
contains only finitely many subgroups, the set $W\!K\backslash
(G/H^K)$ is finite for each $K \subseteq G$ and is non-empty for only
finitely many conjugacy classes $(K)$ of subgroups $K \subseteq G$.
This shows that $\mor(y,x) \not= \emptyset$  for at most finitely many
isomorphism classes $[y]$ of objects $y \in \Pi_0(G,X)$ and that the
$W\!K_y$-set $\mor(y,x)$ decomposes into finitely many  
$W\!K_y$ orbits with finite isotropy groups for each object $y \in
\Pi_0(G,X)$. We define 
\begin{eqnarray}
\ch^G(X)([x])_{[y]} & := & \sum_{\substack{W\!K_y \cdot \sigma \in\\
W\!K_y\backslash \mor(y,x)}}~  |(W\!K_y)_{\sigma}|^{-1},
\label{definition of chi^G(X)([x])_{[y]}}
\end{eqnarray}
where $(W\!K_y)_{\sigma}$ is the isotropy group of $\sigma \in
\mor(y,x)$ under the $W\!K_y$-action.

\begin{lem} 
\label{lem: injectivity of the character map chi^G}
Let $X$ be a finite proper $G$-$CW$-complex. Then the
map $\ch^G(X)$ of \eqref{definition of character map chi^G(X)} is
injective and satisfies 
$$\ch^G(X)(\chi^G(X))_{[y]} ~ = ~ \chi^{\qq W\!K_y}(X^K(y)).$$
The induced map
$$
\id_{\qq} \otimes_{\zz} \ch^G(X) \co \qq \otimes_{\zz} U^G(X)
~ \xrightarrow{\cong} ~
\bigoplus_{(H)\in \consub(G)} ~ \bigoplus_{W\!H\backslash \pi_0(X^H)} \qq$$
is bijective.
\end{lem}
\proof Injectivity of $\chi^G(X)$ and 
$\ch^G(X)(\chi^G(X))_{[y]} ~ = ~ \chi^{\qq W\!K_y}(X^K(y)).$ is proved in
\cite[Lemma 5.3]{Lueck-Rosenberg(2002a)}. The bijectivity of 
$\id_{\qq} \otimes_{\zz} \ch^G(X)$ follows since its source and its
target are  $\qq$-vector spaces of the same finite
$\qq$-dimension. \qed

Now let us briefly summarize the various notions of equivariant Euler
characteristic and the relations among them. Since some of these are
only defined when $M$ is compact and $G$ is finite, we temporarily
make these assumptions \emph{for the rest of this section only}.

\begin{defn}
\label{def: Burnside ring}
If $G$ is a finite group, the \emph{Burnside ring} $A(G)$ of $G$ is the
Grothendieck group of the (additive) monoid of finite $G$-sets, where
the addition comes from disjoint union. This becomes a ring under the
obvious multiplication coming from the Cartesian product of $G$ sets.
There is a natural map of rings $j_1\co A(G)\to R_\bR(G)=KO^G_0(\pt)$ 
that comes from sending a finite $G$-set $X$ to the orthogonal
representation of $G$ on the finite-dimensional real Hilbert space
$L^2_\bR(X)$. This map can fail to be injective or 
fail to be surjective, even
rationally. 
\end{defn}

The rank of $A(G)$ is the number of
conjugacy classes of subgroups of $G$, while the rank of $R_\bR(G)$ is
the number of $\bR$-conjugacy classes in $G$, where $x$ and $y$ are
called $\bR$-conjugate if they are conjugate or if $x^{-1}$ and $y$
are conjugate \cite[\S13.2]{Serre(1977)}. Thus $\rank
A\bigl((\bZ/2)^3\bigr) = 16 > \rank R_\bR\bigl((\bZ/2)^3\bigr) = 8$;
on the other hand, $\rank A(\bZ/5)=2 < \rank R_\bR(\bZ/5)=3$.

\begin{defn}
\label{def: more Euler characteristics}
Now let $G$ be a finite group, $M$ a compact $G$-manifold (without
boundary). We define three more equivariant Euler characteristics for
$G$:
\begin{enumerate}
\item the analytic equivariant Euler characteristic $\chi^G_a(M)
\in KO^G_0$, the
equivariant index of the Euler operator on $M$. Since the index of an
operator is computed by pushing its $K$-homology class forward to
$K$-homology of a point, $\chi^G_a(M)=c_*(\Eul^G(M))$, where $c\co M
\to \pt$ and $c_*$ is the induced map on $KO_0^G$.
\item the stable homotopy-theoretic 
equivariant Euler characteristic $\chi^G_s(M)\in A(G)$. 
This is discussed, say, in
\cite{Dieck(1987)}, Chapter IV, \S2.
\item a
certain unstable homotopy-theoretic equivariant
Euler characteristic, which we will denote here $\chi^G_u(M)$ to distinguish
it from $\chi^G_s(M)$. This invariant is defined in \cite{WW2}, and
shown to be the obstruction to existence of an everywhere
non-vanishing $G$-invariant  vector field on $M$. The invariant
$\chi^G_u(M)$ lives in a group $A^G_u(M)$ (Waner and Wu call it
$A_M(G)$, but the notation $A^G_u(M)$ is more consistent with our
notation for $U^G(M)$) defined as follows: $A^G_u(M)$ is the free
abelian group on finite $G$-sets \emph{embedded in $M$}, modulo
isotopy (if $s_t$ is a $1$-parameter family of finite $G$-sets
embedded in $M$, all isomorphic to one another as $G$-sets, then
$s_0\sim s_1$) and the relation $[s\amalg t]=[s]+[t]$. Waner and Wu 
define a map $d\co A^G_u(M)\to A(G)$ (defined by forgetting that a
$G$-set $s$ is embedded in $M$, and just viewing it abstractly)
which maps $\chi^G_u(M)$ to $\chi^G_s(M)$. Both $\chi^G_u(M)$ 
and $\chi^G_s(M)$ may be computed from the virtual finite $G$-set given
by the singularities of a $G$-invariant canonically transverse vector field,
where the signs are given by the indices at the singularities.
\end{enumerate}
\end{defn}
\begin{prop}
\label{prop: relations among Euler characteristics}
Let $G$ be a finite group, and let $M$ be a compact $G$-manifold 
{\lp}without boundary{\rp}. The following diagram commutes:
\begin{equation}
\label{eq: diagram of Euler characteristics}
\xymatrix{
A_u^G(M) \ar[r] \ar[dr]_d &  U^G(M)  \ar[r]^{e^G(M)}  \ar[d]^{c_*}
& KO^G_0(M) \ar[d]^{c_*}\\
& A(G)=U^G(\pt) \ar[r]^(.45){j_1}  & \,R_\bR(G)=KO^G_0(\pt).}
\end{equation}
The map $A_u^G(M) \to U^G(M)$
in the upper left is an isomorphism if $(M,\,G)$ satisfies the
\emph{weak gap hypothesis}, that is, if whenever 
$H\subsetneq K$ are subgroups
of $G$, each component of $G^K$ has codimension at least $2$ in 
the component of $G^H$ that contains it \textup{\cite{WW2}}. 
Furthermore, under the maps of this diagram,
\[
\begin{aligned}
\chi^G_u(M)&\mapsto \chi^G(M),\\
c_*\co \chi^G(M)&\mapsto \chi^G_s(M),\\
d\co \chi^G_u(M)&\mapsto \chi^G_s(M),
\end{aligned}
\qquad 
\begin{aligned}
e^G(M)\co\chi^G(M)&\mapsto
\Eul^G(M), \\
j_1\co \chi^G_s(M)&\mapsto
\chi^G_a(M),\\
c_*\co \Eul^G(M)&\mapsto \chi^G_a(M).
\end{aligned}
\]
\end{prop}
\begin{proof}
This is just a matter of assembling known information. The facts about
the map $A_u^G(M) \to U^G(M)$ are in \cite[\S2]{WW2} and in
\cite{Lueck-Rosenberg(2002a)}. That $e^G(M)$ sends $\chi^G(M)$ to
$\Eul^G(M)$ is Theorem \ref{the: e^G(M)(chi^G(M)) = Eul^G(M)}.
Commutativity of the square follows immediately from the definition of
$e^G(M)$, since $c_*\circ e^G(M)$ sends the basis element associated
to  $(H) \subseteq \consub(G)$ 
and $W\!H \cdot C \in W\!H\backslash \pi_0(M^H)$ to the class of the
orthogonal representation of $G$ on $L^2(G/H)$. But under $c_*$, this
same basis element maps to the $G$-set $G/H$ in $A(G)$, which 
also maps to the
orthogonal representation of $G$ on $L^2(G/H)$ under $j_1$.
\end{proof}


\tit{The transformation $e^G(X)$}{The transformation eG(M)}

Next we factorize the transformation $e^G(M)$ defined in
\eqref{definition of e^G(M)} as
\begin{multline*}
e^G(M)\co U^G(M) \xrightarrow{e_1^G(M)} H_0^{\Or(G)}(M;\underline{R_{\qq}})
\xrightarrow{e_2^G(M)}
H_0^{\Or(G)}(M;\underline{R_{\rr}})
\\  \xrightarrow{e_3^G(M)}
KO_0^G(M),
\end{multline*}
where $e_2^G(X)$ and  $e_3^G(X)$ are  rationally injective. Rationally
we will identify 
$H_0^{\Or(G)}(M;\underline{R_{\qq}})$ and the element $e^G_1(M)(\chi^G(M))$
in terms of the orbifold Euler characteristics $\chi^{W\!L_C}(C)$,
where $(L)$ runs through the  
conjugacy classes of finite cyclic subgroups $L$ of $G$ and $W\!L \cdot C$ 
runs through the orbits in $W\!L\backslash \pi_0(X^L)$. Here 
$W\!L_C$ is the isotropy group of 
$C \in \pi_0(X^L)$ under the $W\!L=N\!L/L$-action. 
Notice that $e^G_1(M)(\chi^G(M))$ carries rationally the same
information as $\Eul^G(M) \in KO_0^G(M)$. 
 
Here and elsewhere $H_0^{\Or(G)}(M;\underline{R_F})$ is the Bredon
homology of $M$ with coefficients the covariant functor 
$$\underline{R_F}\co \Orcat \to \zz\text{-}\Mod.$$
The orbit category $\Orcat$ 
has as objects homogeneous spaces $G/H$ with finite $H$
and as morphisms $G$-maps (Since $M$ is proper,
it suffices to consider coefficient systems over $\Orcat$ instead over
the full orbit category $\Or(G)$.)
The functor  $\underline{R_F}$ into the category
$\zz\text{-}\Mod$ of $\zz$-modules
sends $G/H$ to the representation ring
$R_F(H)$ of the group $H$ over the field $F = \qq$, $\rr$ or $\cc$. It sends
a morphism $G/H \to G/K$ given by $g'H \mapsto g'gK$ for some
$g \in G$ with $g^{-1}Hg \subseteq K$ to the induction homomorphism
$R_F(H) \to R_F(K)$ associated with the group homomorphism
$H \to K, h \mapsto g^{-1}hg$.  This is independent of the choice of
$g$ since an inner automorphism of $K$ induces the identity on $R_{\rr}(K)$.
Given a covariant functor $V\co  \Or(G) \to \zz\text{-}\Mod$, the
\emph{Bredon homology} of a $G$-$CW$-complex $X$ with coefficients 
in $V$ is defined as follows.
Consider the cellular (contravariant) $\zz\Or(G)$-chain complex 
$C_*(X^-)\co \Or(G) \to \zz\text{-}\Chain $ which assigns to $G/H$ the
cellular chain  complex of the $CW$-complex $X^H = \map_G(G/H,X)$.  
One can form the tensor product over the orbit category
(see for instance \cite[9.12 on page 166]{Lueck(1989)})
$C_*(X^-) \otimes_{\zz\Orcat} V$ which is a $\zz$-chain complex and 
whose homology groups are defined to be $H_p^{\Or(G)}(X;V)$. 

The zero-th Bredon homology can be made more explicit.
Let
\begin{eqnarray}
Q\co \Pi_0(G;X) \to \Orcat
\end{eqnarray}
be the forgetful functor sending an object $x\co G/H \to X$ to $G/H$. 
Any covariant functor $V\co \Orcat \to \zz\text{-}\MOD$ induces a functor
$Q^*V\co \Pi_0(G;X) \to \zz\text{-}\Mod$ by composition with $Q$.
The colimit (= direct limit) of the functor $Q^*\underline{R_F}$ 
is naturally isomorphic to the Bredon homology
\begin{eqnarray}
\beta^G_F(X)\co \colim_{\Pi_0(G,X)} Q^* R_F(H) &  
\xrightarrow{\cong}  &
H_0^{\Or(G)}(X;\underline{R_F}).
\label{identification beta^G_1(X) of colimit and H_0^{Or(G)}}
\end{eqnarray}
The isomorphism $\beta^G_F(X)$ above is induced by the various maps
$$R_F(H) = H_0^{\Or(G)}(G/H;\underline{R_F})
\xrightarrow{H_0^{\Or(G)}(x;\underline{R_F})} 
H_0^{\Or(G)}(X;\underline{R_F}),$$
where $x$ runs through all $G$-maps $x\co G/H \to X$. We define natural maps 
\begin{eqnarray}
e_1^G(X)\co  U^G(X) & \to  & H_0^{\Or(G)}(X;\underline{R_{\qq}});
\label{map e^G_1(X) from U^g(X) to Bredon}
\\
e_2^G(X)\co  H_0^{\Or(G)}(X;\underline{R_{\qq}}) & \to  &
H_0^{\Or(G)}(X;\underline{R_{\rr}}); 
\label{passage e^G_2(X) from qq to rr}
\\
e_3^G(X)\co  H_0^{\Or(G)}(X;\underline{R_{\rr}}) & \to & KO_0^G(X)
\label{edge homomorphism e^G_3(X)}
\end{eqnarray}
as follows. The map $\beta^G_{\qq}(X)^{-1} \circ  e^G_1(X)$  sends the
basis element 
$[x\co  G/H \to X]$ to the image of the trivial representation $[\qq]
\in R_{\qq}(H)$ 
under the canonical map associated to $x$ 
$$R_{\qq}(H) ~ \to ~ \colim_{\Pi_0(G,X)} Q^*R_{\qq}(H).$$ 
The map $e^G_2(X)$ is induced by the change of fields homomorphisms
$R_{\qq}(H) \to R_{\rr}(H)$ for $H \subseteq G$ finite.
The map $e_3^G(X) \circ \beta^G(X)$ is the colimit over the system of maps
$$R_{\rr}(H) = KO_0^H(\{\ast\}) \xrightarrow{(\alpha_H)_*} KO_0^G(G/H)
\xrightarrow{KO_0^G(x)} KO_0^G(X)$$ 
for the various $G$-maps $x\co  G/H \to X$, where $\alpha_H\co  H \to
G$ is the inclusion.

\begin{thm} \label{the: e^G(X)} Let $X$ be a proper $G$-$CW$-complex. Then
\begin{enumerate}

\item \label{the: e^G(X): e^G(X) = e^G_3(X) circ e_2^G(X) circ e_1^G(X)}
The map $e^G(X)$ defined in \eqref{definition of e^G(M)} factorizes as
\begin{multline*}e^G(X)\co  U^G(X) 
 \xrightarrow{e_1^G(X)}  H_0^{\Or(G)}(X;\underline{R_{\qq}})
\xrightarrow{e_2^G(X)}   
H_0^{\Or(G)}(X;\underline{R_{\rr}}) \\
\xrightarrow{e_3^G(X)}   KO_0^G(X);
\end{multline*}

\item \label{the: e^G(X): rational injectivity of e^G_2(X)}
The map 
$$\qq \otimes_{\zz} e^G_2(X) \co \qq \otimes_{\zz}
H_0^{\Or(G)}(X;\underline{R_{\qq}}) \to 
\qq \otimes_{\zz} H_0^{\Or(G)}(X;\underline{R_{\rr}})$$
is injective;

\item \label{the: e^G(X): Chern character}
For each $n \in \zz$ there is an isomorphism, natural in $X$,
$$\chern^G_n(X)\co \bigoplus_{\substack{p,q \in \zz,\\p \ge 0,p+q = n}} 
\qq \otimes_{\zz} H_p^{\Or(G)}(X;\underline{KO}_q^G) \xrightarrow{\cong}
\qq \otimes_{\zz} KO_n^G(X),$$
where $\underline{KO}^G_q$ is the covariant functor from $\Orcat$ to
$\zz\text{-}\Mod$ sending $G/H$ to $KO_q^G(G/H)$. The map
$$\id_{\qq} \otimes_{\zz} e^G_3(X)\co  \qq \otimes_{\zz}
H_0^{\Or(G)}(M;\underline{R_{\rr}})  
\to  \qq \otimes_{\zz} KO_0^G(M)$$
is the restriction of $\chern^G_n(X)$ to the summand for
$p = q = 0$ and is hence injective;

\item \label{the: e^G(X): injectivity of e^G_3(X) for small dim}
Suppose $\dim(X) \le 4$ and one of the following conditions is
satisfied:
\begin{enumerate}
\item either $\dim(X^H) \le 2$ for $H\subseteq G$, 
$H \not= \{1\}$, or
\item no subgroup $H$ of $G$ has irreducible representations of
complex or quaternionic type.
\end{enumerate}
Then 
$$e^G_3(X)\co  H_0^{\Or(G)}(M;\underline{R_{\rr}}) 
\to   KO_0^G(M)$$
is injective.

\end{enumerate}
\end{thm}

\proof 
\eqref{the: e^G(X): e^G(X) = e^G_3(X) circ e_2^G(X) circ e_1^G(X)}
follows directly from the
definitions. 

\eqref{the: e^G(X): rational injectivity of e^G_2(X)} will be proved later.

\eqref{the: e^G(X): Chern character}\qua
An equivariant Chern character $\chern^G_*$ for equivariant homology
theories such as 
equivariant K-homology is constructed in 
\cite[Theorem 0.1]{Lueck(2002b)}
(see also \cite[Theorem 0.7]{Lueck(2002c)}). The restriction of
$\chern^G_0(M)$ to $H_0^{\Or(G)}(X;\underline{R_{\rr}})$  is just
$e^G_3(X)$ under the identification 
$\underline{R_{\rr}} = \underline{KO}^G_0$.

\eqref{the: e^G(X): injectivity of e^G_3(X) for small dim}\qua
Consider the equivariant
Atiyah-Hirzebruch spectral sequence which converges to $KO_{p+q}^G(X)$
(see for instance \cite[Theorem 4.7 (1)]{Davis-Lueck(1998)}). Its 
$E^2$-term is $E^2_{p,q} = H_p^{\Or(G)}(X;\underline{KO_q^G})$.
The abelian group $KO_q^G(G/H)$ is isomorphic to the real topological
$K$-theory  
$KO_q(\rr H)$ of the real $C^*$-algebra $\rr H$. The  real $C^*$-algebra
$\rr H$ splits as a product of matrix algebras over $\rr$, $\cc$ and
$\hh$, with as many summands of a given type as there are irreducible
real representations of $H$ of that type (see \cite[\S13.2]{Serre(1977)}).
By Morita invariance of topological real $K$-theory, we conclude that
$KO_q^G(G/H)$ is a direct sum of copies of the non-equivariant K-homologies
$KO_q(\ast) = KO_q(\rr)$, $KU_q(\ast) = KO_q(\cc)$, and
$KSp_q(\ast) = KO_q(\hh)$. In particular, we conclude 
that $KO_q^G(G/H) = 0$ for $q\equiv -1$ (mod $8$).
As a consequence, $\underline{KO_{-1}^G}=0$ and
$E^2_{p,-1}=H_p^{\Or(G)}(X;\underline{KO_{-1}^G})=0$.
If no subgroup $H$ of $G$ has irreducible representations of
complex or quaternionic type, then similarly
$KO_q^G(G/H) = 0$ for all subgroups $H$ of $G$ and $q=-2,\,-3$,
and $E^2_{p,-2}=0$, $E^2_{p,-3}=0$, as well.

Let $X^{>1}\subseteq X$ be the subset $\{x \in X \mid G_x \not=1\}$.
There is a short exact sequence
$$H_p^{\Or(G)}(X^{>1};\underline{KO_q^G}) \to
H_p^{\Or(G)}(X;\underline{KO_q^G}) \to 
 H_p^{\Or(G)}(X,X^{>1};\underline{KO_q^G}).$$
Since the isotropy group of any point in $X-X^{>1}$ is trivial, we get
an isomorphism 
$$H_p^{\Or(G)}(X,X^{>1};\underline{KO_q^G}) = H_p(C_*(X,X^{>1})
\otimes_{\zz G} KO_q^G(G/1))$$ 
Since $KO_q^G(G/1)) = KO_q(\rr)$ vanishes for $q \in \{-2,-3\}$, we get
for $p \in \zz$ and  $q \in \{-2,-3\}$
$$H_p^{\Or(G)}(X,X^{>1};\underline{KO_q^G}) = 0.$$
So if
$\dim(X^H) \le 2$ for $H\subseteq G$, $H \not= \{1\}$, we have
$\dim(X^{>1}) \le 2$.  This implies for $p \ge 3$ and $q \in \zz$
$$H_p^{\Or(G)}(X^{>1};\underline{KO_q^G}) = 0.$$
We conclude that for $p \ge 3$ and $q \in \{-2,-3\}$
$$E^2_{p,q} = H_p^{\Or(G)}(X;\underline{KO_q^G}) = 0,$$
just as in the previous case.

Now there are no non-trivial differentials out of
$E^r_{0,0}$, and since $\dim X\le 4$,
$d^r_{r,1-r}\co E^r_{r,1-r} \to E^r_{0,0}$
must be zero for $r>4$. But we have just seen that
$E^2_{2,-1}=0$, $E^2_{3,-2}=0$, and $E^2_{4,-3}=0$. 
Hence each differential $d^r_{p,q}$ 
which has $E^r_{0,0}$ as source or target is trivial. 
Hence the edge homomorphism 
restricted to $E^2_{0,0}$ is injective. But this map is  $e^G_3(X)$.
\qed

\begin{rem} \label{rem: comparing chi^G(X) and Eul^G(M)}
We conclude from Theorem \ref{the: e^G(M)(chi^G(M)) = Eul^G(M)}
and Theorem  \ref{the: e^G(X)} that $\Eul^G(M)$ carries rationally the same
information as the image of the equivariant Euler characteristic
$\chi^G(X)$ under the map 
$e^G_1(X)\co U^G(X) \to H_0^{\Or(G)}(M;\underline{R_{\rr}})$.
Moreover, in contrast to the class of  the signature operator, the class
$\Eul^G(M) \in KO_0^G(M)$ of the Euler operator does not carry ``higher''
information 
because its preimage under the equivariant Chern character is concentrated 
in the summand corresponding to $p = q = 0$.
\end{rem}

Next we recall the definition of the Hattori-Stallings 
rank of a finitely generated projective $RG$-module $P$,
for some commutative ring $R$ and a group $G$. Let $R[\con(G)]$ be the
$R$-module with 
the set of conjugacy classes $\con(G)$ of elements in $G$ as basis.
Define the universal $RG$-trace
$$\tr_{RG}^u\co RG \to R[\con(G)], \hspace{5mm} \sum_{g \in G} r_g
\cdot g ~ \mapsto ~  
\sum_{g \in G} r_g \cdot (g).$$
Choose a matrix $A = (a_{i,j}) \in M_n(RG)$
such that $A^2 = A$ and the image of the map $r_A\co RG^n \to RG^n$
sending $x$ to $xA$ is 
$RG$-isomorphic to $P$. Define the  \emph{Hattori-Stallings rank} 
\begin{eqnarray}
\HS_{RG}(P) & := \sum_{i=1}^n \tr_{RG}^u(a_{i,i}) \hspace{5mm} \in
R[\con(G)]. 
\label{def of Hattori-Stallings rank for R_{rr}(H)}
\end{eqnarray}
Let $\alpha\co H_1 \to  H_2$ be a group homomorphism. It induces a map
\[
\con(H_1) \to \con(H_2), \quad (h) \mapsto (\alpha(h))
\]
and thus an $R$-linear map $\alpha_*\co R[\con(H_1)] \to R[\con(H_2)]$. 
If $\alpha_*P$ is the $R[H_2]$-module 
obtained by induction from the finitely generated projective
$R[H_1]$-module $P$, then  
\begin{eqnarray}
\HS_{RH_2}(\alpha_*P) & = & \alpha_*(\HS_{RH_1}(P)).
\label{Hattori-Stallings rank and induction}
\end{eqnarray}

Next we compute 
$F \otimes_{\zz} \left(\colim_{\Pi_0(G,X)} Q^*R_{F}(H)\right)$ for $F
= \qq, \rr, \cc$.
Two elements $g_1$ and $g_2$ of a group $G$ are called \emph{$\qq$-conjugate}
if and only if the cyclic  subgroup $\langle g_1\rangle $ generated by
$g_1$ and  
the cyclic  subgroup $\langle g_2 \rangle$ generated by $g_2$ are
conjugate in $G$. 
Two elements $g_1$ and $g_2$ of a group $G$ are called \emph{$\rr$-conjugate}
if and only if $g_1$ and $g_2$  or $g_1^{-1}$ and $g_2$ are conjugate in $G$.
Two elements $g_1$ and $g_2$ of a group $G$ are called \emph{$\cc$-conjugate}
if and only if they are conjugate in $G$ (in the usual sense).
We denote by $(g)_{F}$ the set of elements of $G$ which are
$F$-conjugate to $g$. 
Denote by $\con_{F}(G)$ the set of $F$-conjugacy classes $(g)_F$ of
elements of finite 
order $g \in G$. Let $\class_{F}(G)$ be the $F$-vector space generated
by the set $\con_{F}(G)$. 
This is the same as the $F$-vector space of functions $\con_{F}(G)  \to \rr$ 
whose support is finite.

Let $\pr_F\co \con(H) \to \con_F(H)$ be the canonical
epimorphism for a finite group $H$. It extends to an
$F$-linear epimorphism $F[\pr_F]\co F[\con(H)] \to \class_F(H)$.
Define for a finite-dimensional $H$-representation $V$ over $F$ for a finite group $H$
\begin{eqnarray}
\HS_{F,H}(V) & := & F[\pr_F](\HS_{FH}(V))\in  \class_F(H).
\label{definition of Hs_{F,H}(V)}
\end{eqnarray}
Let $\alpha\co H_1 \to  H_2$ be a homomorphism of finite groups.
It induces a map $\con_F(H_1) \to \con_F(H_2), (h)_F \mapsto
(\alpha(h))_F$ and thus an $F$-linear map 
\[
\alpha_* \co \class_F(H_1) \to \class_F(H_2). 
\]
If $V$ is a finite-dimensional $H_1$-representation over $F$, we
conclude from 
\eqref{Hattori-Stallings rank and induction}:
 \begin{eqnarray}
\HS_{F,H_2}(\alpha_*V) & = & \alpha_*(\HS_{F,H_1}(V)).
\label{Hattori-Stallings rank and induction for representations over F}
\end{eqnarray}

\begin{lem} \label{R_F(H) and class_F(H)} 
Let $H$ be a finite group and $F = \qq, \rr$ or $\cc$. 
Then the Hattori-Stallings rank defines an isomorphism
$$\HS_{F,H}\co F \otimes_{\zz} R_F(H) ~ \xrightarrow{\cong} \class_F(H)$$
which is natural with respect to induction with respect to group homomorphism
$\alpha\co H_1 \to H_2$ of finite groups.
\end{lem}
\proof 
One easily checks for a finite-dimensional $H$-representation over $F$
of a finite group $H$ 
$$\HS_{FH}(V) ~ = ~ \sum_{(h) \in \con(H)} \frac{|(h)|}{|H|} \cdot
\tr_F(l_h).$$ 
This explains the relation between the Hattori-Stallings rank and the
character of a representation --- they contain equivalent information. 
(We prefer the Hattori-Stallings rank because it behaves better under
induction.) 
We conclude from \cite[page 96]{Serre(1977)} that $\HS_{FH}(V)$ as a function
$\con(H) \to F$ is constant on the $F$-conjugacy classes of elements in $H$ and 
that  $\HS_{F,H}$ is bijective.
Naturality follows from \eqref{Hattori-Stallings rank and induction for representations over F}.
\qed

Let $\underline{\class_F}$ be  the covariant functor $\Orcat \to F\text{-}\MOD$
which sends an object $G/H$ to $\class_F(H)$.
The isomorphisms $\HS_{F,H}\co  F \otimes_{\zz} R_F(H) ~ \xrightarrow{\cong} \class_F(H)$ 
yield a natural equivalence of covariant functors from $\pi_0(G,X)$ to $F\text{-}\Mod$.
Thus we obtain an isomorphism
\begin{multline}
\HS^G_F(X)\co 
F \otimes_{\zz} \colim_{\Pi_0(G,X)}  Q^*\underline{R_{F}}
\\
\xrightarrow{\cong} \colim_{\Pi_0(G,X)}  Q^*F \otimes_{\zz} \underline{R_{F}}
\xrightarrow{\cong} \colim_{\Pi_0(G,X)}  Q^*\underline{\class_{F}}.
\label{isomorphism HS^G(X) induced by Hattori-Stallings rank}
\end{multline}

Let $f\co S_0 \to  S_1$ be a map of sets. It extends to an $F$-linear map
$F[f]\co F[S_0] \to F[S_1]$. Suppose that the  
preimage of any element in $S_1$ is finite. Then we obtain an $F$-linear map
\begin{eqnarray}
f^*\co  F[S_1] & \to & F[S_0], \hspace{5mm} 
s_1 \mapsto \sum_{s_0 \in f^{-1}(s_1)} s_0.
\label{definition of f^* for finite-one maps of sets}
\end{eqnarray}
If we view elements in $F[S_i]$ as functions $S_i \to F$, then $f^*$ is
given by composing 
with $f$. One easily checks that $F[f] \circ f^*$ is bijective and
that for a second map 
$g\co  S_1 \to S_2$, for which the preimages of any element in $S_2$
is finite, we have $f^* \circ g^* = (g \circ f)^*$. 

Now we can finish  the proof of Theorem \ref{the: e^G(X)} by explaining how assertion 
\eqref{the: e^G(X): rational injectivity of e^G_2(X)} is proved.\\
\proof 
Let $H$ be a finite group. Let $p_H\co \con_{\rr}(H) \to \con_{\qq}(H)$
be the projection. If $V$ is a finite-dimensional 
$H$-representation over $\qq$, then $\rr \otimes_{\qq} V$ is a
finite-dimensional 
$H$-representation over $\rr$ and $\HS^{\rr H}(\rr \otimes_{\qq} V)$
is the image of 
$\HS^{\qq H}(V)$ under the obvious map
$\qq[\con(H)] \to \rr[\con(H)]$. Recall that $\HS_{F,H}(V)$ is the
image of $\HS^{FH}(V)$ under 
$F[\pr_F]$ for $\pr_F \co \con(H) \to \con_F(H)$ the canonical projection; 
$\HS^{FH}(V)$ is constant on the $F$-conjugacy classes of elements in $H$.
This implies that the following diagram commutes
$$
\begin{CD}
\rr \otimes_{\zz} R_{\qq}(H) @ > \id_{\rr}  \otimes_{\zz}\HS_{\qq,H} > \cong >  
\rr \otimes_{\qq} \class_{\qq}(H)
\\
@V \ind(H) VV @V \ind(H) VV 
\\
\rr \otimes_{\zz} R_{\rr}(H) @ > \HS_{\rr,H} > \cong >  \,\class_{\rr}(H),
\end{CD}
$$
where the left vertical arrow comes from inducing a
$\qq$-representation to a $\rr$-representation and 
the right vertical arrow is 
$$\rr \otimes_{\qq} \class_{\qq}(H) =
\rr \otimes_{\qq} \qq[\con_{\qq}(H)] \xrightarrow{\id_{\rr}
\otimes_{\qq} (p_H)^*} 
\rr \otimes_{\qq} \qq[\con_{\rr}(H)] = \class_{\rr}(H).$$
Let
$$q(H)\co  \class_{\rr}(H) \to \rr \otimes_{\qq} \class_{\qq}(H)$$
be the map 
$$\rr[p_H]\co \rr[\con_{\rr}(H)]\to \rr[\con_{\qq}(H)] = \rr
\otimes_{\qq} \qq[\con_{\qq}(H)]$$ 
The $q(H) \circ \ind(H)$ is bijective.

We get natural transformations of functors from $\Orcat$ to
$\rr\text{-}\Mod$: 
\begin{eqnarray*} 
\ind \co \underline{\rr \otimes_{\zz} R_{\qq}}   & \to & 
\underline{\rr \otimes_{\zz} R_{\rr}},
\\
\ind \co \underline{\rr \otimes_{\qq} \class_{\qq}}   & \to &
\underline{\class_{\rr}},
\\
q \co \underline{\class_{\rr}} & \to & \underline{\rr \otimes_{\qq}
\class_{\qq}}. 
\end{eqnarray*}
Also $q \circ \ind$ is a natural equivalence. 
This implies that $\ind$ induces a split injection on the colimits
$$\colim_{\Pi_0(G;X)} Q^*\ind \co
\colim_{\Pi_0(G;X)} Q^*\underline{\rr \otimes_{\qq} \class_{\qq}} \to
\colim_{\Pi_0(G;X)} Q^*\underline{\class_{\rr}}.$$
Since the following diagram commutes and has isomorphisms as vertical arrows
$$
\begin{CD}
\rr \otimes_{\zz} H_0^{\Or(G)}(X;\underline{R_{\qq}})
@>\rr \otimes_{\zz}  e^G_2(X) >>
\rr \otimes_{\zz} H_0^{\Or(G)}(X;\underline{R_{\rr}})
\\
@A \rr \otimes_{\zz} \beta^G_{\qq}(X) A\cong A 
@A \rr \otimes_{\zz} \beta^G_{\rr}(X) A\cong A 
\\
\colim_{\Pi_0(G;X)} Q^*\underline{\rr \otimes_{\zz} R_{\qq}} 
@> \colim_{\Pi_0(G;X)} Q^*\ind >>
\colim_{\Pi_0(G;X)} Q^*\underline{\rr \otimes_{\zz} R_{\rr}}
\\
@V \rr \otimes_{\qq} \HS^G_{\qq}(X) V \cong V
@V \rr \otimes_{\zz} \HS^G_{\rr}(X) V \cong V
\\
\colim_{\Pi_0(G;X)} Q^*\underline{\rr \otimes_{\qq} \class_{\qq}} 
@> \colim_{\Pi_0(G;X)} Q^*\ind >>
\,\,\colim_{\Pi_0(G;X)} Q^*\underline{\class_{\rr}}\,,
\end{CD}
$$
the top horizontal arrow is split injective. Hence $e^G_2(X)$ is
rationally split injective. 
This finishes the proof of Theorem \ref{the: e^G(X)}
\eqref{the: e^G(X): rational injectivity of e^G_2(X)}. \qed

\begin{notation} \label{not: ZH and ZH_y}
Consider $g \in G$ of finite order. Denote by 
$\langle g \rangle$ the finite cyclic subgroup generated by $g$. 
Let $y\co G/\langle g \rangle \to X$ be a $G$-map.
Let $F$ be one of the fields $\qq$, $\rr$ and $\cc$. Define
\begin{eqnarray*}
C_{\qq}(g) & = & \{g' \in G, (g')^{-1}gg' \in \langle g \rangle\}\};
\\
C_{\rr}(g) & = & \{g' \in G, (g')^{-1}gg' \in \{g,g^{-1}\}\};
\\
C_{\cc}(g) & = & \{g' \in G, (g')^{-1}gg' = g\}.
\end{eqnarray*}
Since $C_F(g)$
is a subgroup of the normalizer $N\!\langle g \rangle$ of $\langle g \rangle$ in $G$
and contains $\langle g \rangle$, we can define a  subgroup $Z_F(g) \subseteq W\!\langle g \rangle$ by
$$Z_F(g)  ~ := ~ C_F(g)/\langle g \rangle.$$
Let $Z_F(g)_y$ be the intersection of $W\!\langle g \rangle_y$ (see
Notation \ref{not: X^H(x) and so on}) with $Z_F(g)$, or, equivalently,
the subgroup of $Z_F(g)$ represented by elements $g \in  C_F(g)$ for which 
$g\cdot y(1 \langle g \rangle)$ and $y(1 \langle g \rangle)$ 
lie in the same component of $X^{\langle g \rangle}$.
\end{notation}

It is useful to interpret the group $C_F(g)$ as follows. Let $G$ act on the set
$G_f := \{g\mid g \in G, |g| < \infty\}$ by conjugation. Then the projection
$G_f \to \con_{\cc}(G)$ induces a bijection
$G\backslash G_f \xrightarrow{\cong} \con_{\cc}(G)$ and the isotropy
group of $g \in G_f$ is $C_{\cc}(g)$. Let
$G_f/\{\pm 1\}$ be the quotient of $G_f$ under the $\{\pm 1\}$-action given by $g \mapsto g^{-1}$.
The conjugation action of $G$ on $G_f$ induces an  action on $G_f/\{\pm 1\}$.
The projection $G \to \con_{\rr}(G) $ induces a bijection
$G\backslash(G/\{\pm 1\})\xrightarrow{\cong} \con_{\rr}(G)$ 
and the isotropy group of $g\cdot\{\pm 1\} \in G_f/\{\pm 1\}$
is $C_{\rr}(g)$. The set $\con_{\qq}(G)$ is the same as the set
$\{(L) \in \consub(G) \mid L \text{ finite cyclic}\}$. For $g \in G_f$ the group
$C_{\qq}(g)$ agrees with $N\langle g \rangle$ and is the isotropy group of
$\langle g \rangle$ in $\{L \subseteq G \mid L \text{ finite cyclic}\}$ under the conjugation action
of $G$. The quotient of $\{L \subseteq G \mid L \text{ finite cyclic}\}$ under the conjugation action
of $G$ is by definition $\con_{\qq}(G) = \{(L) \in \consub(G) \mid L \text{ finite cyclic}\}$.

Consider $(g)_F \in \con_F(G)$. For the sequel we fix a representative $g \in (g)_F$. 
Consider an object of $\Pi_0(G,X)$ of the special form $y\co
G/\langle g \rangle \to  X$. Let
$x\co G/H \to X$ be any object of $\Pi_0(G,X)$.
Recall that $W\!\langle g \rangle_y$ and thus the subgroup $Z_F(g)_y$ act on
$\mor(y,x)$. Define
$$\alpha_F(y,x)\co Z_F(g)_y\backslash \mor(y,x) \to \con_F(H)$$
by sending sending $Z_F(g)_y \cdot \sigma$ for a morphism
$\sigma\co y \to x$, which  given by a $G$-map $\sigma\co G/\langle g
\rangle \to G/H$, 
to $(\sigma(1\langle g \rangle)^{-1} g \sigma(1\langle g \rangle))_F$. 
We obtain a map of sets
\begin{multline}
\alpha_F(x) = \coprod_{(g)_F \in \con(G)_{\rr}} ~
\coprod_{\substack{Z_F(g) \cdot C \in \\
Z_F(g)\backslash\pi_0(X^{\langle g \rangle})}} a(y(C),x) \co  
\\
\coprod_{(g)_F \in \con(G)_F} ~
\coprod_{\substack{Z_F(g)\cdot C \in
\\Z_F(g)\backslash\pi_0(X^{\langle g \rangle})}} 
Z_F(g)_{y(C)}\backslash \mor(y(C),x) \xrightarrow{\cong}  \con_F(H),
\label{definition of alpha_F(x)}
\end{multline}
where we fix for each $Z_F(g) \cdot C \in Z_F(g)\backslash
\pi_0(X^{\langle g \rangle})$  
a representative $C \in \pi_0(X^{\langle g  \rangle})$ and 
$y(C)$ is a fixed morphism $y(C)\co G/\langle g \rangle \to X$ such that
$X^{\langle g \rangle}(y) = C$ in $\pi_0(X^{\langle g \rangle})$. The
map $\alpha_F(x)$ is bijective by the following argument.

Consider $(h)_F \in \con_F(H)$. Let $g \in G$ be the representative of
the class 
$(g)_F$ for which $(g)_F = (h)_F$ holds in $\con_F(G)$. Choose $g_0
\in G$ with $g_0^{-1}gg_0 \in H$ and 
$(g_0^{-1}gg_0)_F = (h)_F$ in $\con_F(H)$. We get a $G$-map
$R_{g_0}\co G/\langle g \rangle \to G/H$ by mapping $g'\langle g
\rangle$ to $g'g_0H$. 
Let $y = y(C) \co G/\langle g \rangle \to X$ be the object chosen above 
for the fixed representative $C$ of $Z_F(g) \cdot X^{\langle g
\rangle}(x \circ R_{g_0})  
\in Z_F(g)_F\backslash \pi_0(X^{\langle g \rangle})$.
Choose $g_1 \in Z_F(g)$ such that the $G$-map
$R_{g_1}\co G/\langle g \rangle \to G/\langle g \rangle$ sending
$g'\langle g \rangle$ to 
$g'g_1\langle g \rangle$ defines a morphism $R_{g_1}\co y \to x \circ
R_{g_0}$ in 
$\Pi_0(G,X)$. Then $\sigma := R_{g_0} \circ R_{g_1}\co y \to x$ is a
morphism such that  
$$(\sigma(1 \langle g \rangle)^{-1} g \sigma(1 \langle g \rangle)_F = (h)_F$$
holds in $\con_F(H)$. This shows
that $\alpha(x)$ is surjective.

Consider for $i = 0,1$ elements $(g_i) \in \con_F(G)$, 
$Z_F(g_i) \cdot C_i \in Z_F(g_i)\backslash\pi_0(X^{\langle g_i \rangle})$  and
$Z_F(g_i)_{y_i} \cdot \sigma_i \in
Z_F(g_i)_{y_i}\backslash \mor(y_i,x)$ for $y_i = y(C_i)$ such that
$$(\sigma_0(1\langle g_0 \rangle)^{-1} g_0 \sigma_0(1\langle g_0 \rangle))_F = 
(\sigma_1(1\langle g_1 \rangle)^{-1} g_1 \sigma_1(1\langle g_1 \rangle))_F$$
 holds in $\con_F(H)$. So we get two elements in the source of $\alpha_F(x)$ which are mapped
to the same element under $\alpha_F(x)$. We have to show that these elements in the source agree.
Choose $g_i' \in G$ such that $\sigma_i$ is given by sending
$g'' \langle g_i \rangle$ to $g''g_i'H$. Then $((g_0')^{-1} g_0 g_0')_F$ and
$((g_1')^{-1} g_1 g_1')_F$ agree in $\con_F(H)$. This implies
$(g_0)_F = (g_1)_F$ in  $\con(G)_F$ and hence $g_0 = g_1$.
In the sequel we write $g = g_0 = g_1$.
Since $((g_0')^{-1} g g_0')_F$ and
$((g_1')^{-1} g g_1')_F$ agree in $\con_F(H)$, there exists $h \in H$ with
$$h^{-1} (g_0')^{-1} g g_0' h \left\{\begin{array}{lll}
\in & \langle (g_1')^{-1} g g_1' \rangle & \text{ if } F = \qq;
\\
\in & \{(g_1')^{-1} g g_1', (g_1')^{-1} g^{-1} g_1'\} & \text{ if } F = \rr;
\\
 =  & (g_1')^{-1} g g_1' & \text{ if } F = \cc.
\end{array}
\right.$$
We can assume without loss of generality that $h = 1$, otherwise replace
$g_0'$ by $g_0'h$. Put $g_2:= g_0'(g_1')^{-1}$. Then $g_2$  is an
element in $Z_F(g)$. 
Let $\sigma_2\co G/\langle  g \rangle \to G/\langle g \rangle$ be the
$G$-map which sends 
$g'' \langle g \rangle$ to $g''g_2\langle g \rangle$. We get the equality of $G$-maps
$\sigma_0 = \sigma_1 \circ \sigma_2$. Since $\sigma_i$ is a morphism
$y_i \to x$ for $i = 0,1$,  we conclude 
$g_i' \cdot x(1H) \in X^{\langle g \rangle}(y_i)$ for $i = 0,1$. This implies
that $g_2 \cdot X^{\langle g \rangle}(y_1) = X^{\langle g \rangle}(y_0)$ in 
$\pi_0(X^{\langle g \rangle})$. This shows
$Z_F(g) \cdot X^{\langle g \rangle}(y_0) = Z_F(g)\cdot X^{\langle g \rangle}(y_1)$
in $Z_F(g)\backslash \pi_0(X^{\langle g  \rangle})$ and hence
$y_0 = y_1$. Write in the sequel $y = y_0 = y_1$. 
The $G$-map $\sigma_2$ defines a morphism $\sigma_2\co y \to y$.
We obtain an equality $\sigma_0 = \sigma_1 \circ \sigma_2$ of morphisms
$y \to x$. We conclude
$Z_F(g)_y \cdot \sigma_0 = Z_F(g)_y\cdot \sigma_1$ 
in $Z_F(g)_y\backslash \pi_0(X^{\langle g \rangle})$. 
This finishes the proof that $\alpha(x)$ is bijective.

Let $\underline{\con_F}$ be the covariant functor from $\Orcat$ to the category of 
finite sets which sends an object $G/H$ to $\con_F(H)$.
The map $\alpha_F(x)$ is natural in $x\co G/H \to X$, in other words,
we get a natural equivalence of functors from $\Pi_0(G,X)$ to the category of finite sets. 
We obtain a bijection of sets
\begin{multline*}
\alpha_F\co  \colim_{x\co G/H \to X  \in \Pi_0(G,X)}
\coprod_{(g)_F \in \con_F(G)} ~
\coprod_{\substack{Z_F\langle g \rangle C \in \\Z_F\langle g \rangle
\backslash\pi_0(X^{\langle g \rangle})}}  
Z_F\langle g \rangle_y\backslash \mor(y,x) \\
\xrightarrow{\cong}
\colim_{\Pi_0(G,X)} Q^*\underline{\con_F}.
\end{multline*}
One easily checks that 
$ \colim_{x\co G/H \to X  \in \Pi_0(G,X)} ~ Z_F\langle g
\rangle_y\backslash \mor(y,x) $  
consists of one element, namely, the one represented by $Z_F\langle g
\rangle_y\cdot 
\id_y \in Z_F\langle g \rangle_y\backslash \mor(y,y)$. Thus we obtain
a bijection 
$$
\alpha^G_F(X)\co \coprod_{(g)_F \in \con_F(G)} ~ 
Z_F\langle g \rangle \backslash\pi_0(X^{\langle g \rangle})
\xrightarrow{\cong} \colim_{\Pi_0(G,X)} Q^*\underline{\con_F},
$$
which sends an element 
$Z_F\langle g \rangle \cdot  C$ in $Z_F\langle g \rangle
\backslash\pi_0(X^{\langle g \rangle})$ 
to the class in the colimit represented by $Z_F\langle g \rangle_y
\cdot \id_y$ in 
$Z_F\langle g \rangle_y \backslash \mor(y,y)$ for any  object 
$y\co  G/\langle g \rangle \to X$ for which $Z_F(g) \cdot X^{\langle g
\rangle}(y) = 
Z_F(g) \cdot C$ holds in $Z_F(g)\backslash \pi_0(X^{\langle g \rangle})$.
It yields an isomorphism of $F$-vector spaces denoted in the same way
\begin{eqnarray} 
\hspace{-5mm} \alpha^G_F(X)\co \bigoplus _{(g)_F \in \con_F(G)} ~ 
\bigoplus_{Z_F(g)\backslash\pi_0(X^{\langle g \rangle})} F
& \xrightarrow{\cong} &
\colim_{\Pi_0(G,X)} Q^*\underline{\class_F}.
\label{colim_{xcolon G/H to X  in Pi_0(G,X)} class_F(H)}
\end{eqnarray}
Let us consider in particular the case $F = \qq$. Recall that $\consub(G)$
is the set of conjugacy classes $(H)$ of subgroups of $G$. Then
$\con_{\qq}(G)$  
is the same as the set $\{(L) \in \consub(G)\mid L \text{ finite
cyclic}\}$ and 
$Z_{\qq}(g)$ agrees with $W\!\langle g \rangle$. Thus 
\eqref{colim_{xcolon G/H to X  in Pi_0(G,X)} class_F(H)} becomes an
isomorphism of $\qq$-vector spaces
\begin{eqnarray}
\alpha^G_{\qq}(X)\co \bigoplus _{\substack{(L) \in \consub(G)\\
L \text{ finite cyclic}}} ~ 
\bigoplus_{W\!L\backslash\pi_0(X^L)} \qq
& \xrightarrow{\cong} &
\colim_{\Pi_0(G,X)} Q^*\underline{\class_{\qq}}.
\label{colim_{xcolon G/H to X  in Pi_0(G,X)} class_{qq}(H)}
\end{eqnarray}
Denote by 
$$\pr\co \bigoplus_{(H)\in \consub(G)} ~ \bigoplus_{W\!H\backslash
\pi_0(X^H)} \qq  
~ \to ~
\bigoplus_{\substack{(L)\in \consub(G)\\ L \text{ finite cyclic}}} ~
\bigoplus_{W\!L\backslash \pi_0(X^L)} \qq $$ 
the obvious projection. Let 
\begin{eqnarray}
\hspace{-5mm} D^G(X) \co 
\bigoplus_{\substack{(L)\in \consub(G)\\ L \text{ finite cyclic}}} ~
\bigoplus_{W\!L\backslash \pi_0(X^L)} \qq 
& \to  &
\bigoplus_{\substack{(L)\in \consub(G)\\ L \text{ finite cyclic}}} ~
\bigoplus_{W\!L\backslash \pi_0(X^L)} \qq
\label{definition of D^G(X)}
\end{eqnarray}
be the automorphism 
$\bigoplus_{\substack{(L)\in \consub(G)\\ L \text{ finite cyclic}}}
\frac{|\Gen(L)|}{|L|} \cdot \id$, 
where $\Gen(L)$ is the set of generators of $L$.

Recall the isomorphisms  $\ch^G(X)$,   
$\beta^G_{\qq}(X)$, $\HS^G_{\qq}(X)^{-1}$, $\alpha^G_{\qq}(X)$ and
$D^G(X)$ from 
\eqref{definition of character map chi^G(X)},
\eqref{identification beta^G_1(X) of colimit and H_0^{Or(G)}},
\eqref{isomorphism HS^G(X) induced by Hattori-Stallings rank}
\eqref{colim_{xcolon G/H to X  in Pi_0(G,X)} class_{qq}(H)}
and \eqref{definition of D^G(X)}. We define
\begin{multline}
\gamma^G_{\qq}(X) \co 
\bigoplus_{\substack{(L)\in \consub(G)\\ L \text{ finite cyclic}}} ~ 
\bigoplus_{W\!L\backslash \pi_0(X^L)} \qq 
~ \xrightarrow{\cong} ~ 
\qq \otimes_{\zz} H_0^{\Or(G)}(X;\underline{R_{\qq}})
\label{definition of gamma^G_{qq}}
\end{multline}
to be the composition $\gamma^G_{\qq}(X) := 
\beta^G_{\qq}(X) \circ \HS^G_{\qq}(X)^{-1} \circ \alpha^G_{\qq}(X)
\circ D^G(X)$. 

\begin{thm} \label{the: commutative diagram for e^G_1(M)}
 The following diagram commutes
\begin{center}
$\begin{CD}
\qq \otimes_{\zz} U^G(X) 
@> \id_{\qq} \otimes_{\zz} e^G_1(X)>> 
\qq \otimes_{\zz} H_0^{\Or(G)}(X;\underline{R_{\qq}})
\\
@V\id_{\qq} \otimes_{\zz} \ch^G(X) V \cong V  
@A\gamma^G_{\qq}(X) A \cong A 
\\
\bigoplus_{(H)\in \consub(G)} ~ \bigoplus_{W\!H\backslash \pi_0(X^H)} \qq  
@> \pr >>
\bigoplus_{\substack{(L)\in \consub(G)\\ L \text{ finite cyclic}}} ~ 
\bigoplus_{W\!L\backslash \pi_0(X^L)} \qq
\end{CD}$
\end{center}
and has isomorphisms as vertical arrows.

The element $\id_{\qq} \otimes_{\zz} e^G_1(X)(\chi^G(X))$ agrees with
the image under the isomorphism $\gamma^G_{\qq}$
of the element 
$$\{\chi^{\qq W\!L_C}(C) \mid (L) \in \consub(G), L \text{ finite cyclic}, 
W\!L\cdot C \in W\!L\backslash \pi_0(X^L)\}$$
given by the various orbifold Euler characteristics of the
$W\!L_C$-$CW$-complexes $C$, 
where $W\!L_C$ is the isotropy group of $C \in \pi_0(X^L)$ under the
$W\!L$-action. 
\end{thm}
\proof 
It suffices to prove the commutativity of the diagram above, then the
rest follows from 
Lemma \ref{lem: injectivity of the character map chi^G}.

Recall that $U^G(X)$ is the free abelian group generated by the set
of isomorphism classes $[x]$ of objects $x\co G/H \to X$.  Hence it
suffices to prove  
for any $G$-map $x\co G/H \to X$
\begin{multline}
\left(\alpha_{\qq}^G(X) \circ D^G(X) \circ \pr \circ \ch^G(X)\right)([x])
\\ ~  = ~ 
\HS^G(X) \circ \beta^G(X)^{-1} \circ e^G_1(X)([x]).
\label{the: commutative diagram for e^G_1(M): 1}
\end{multline}
Given a finite cyclic  subgroup $L \subseteq G$ and a component $C \in \pi_0(X^L)$ 
the element $\left(D^G(X) \circ \pr \circ \ch^G(X)\right)([x])$ has 
as entry in the summand belonging to $(L)$ and $W\!L\cdot C \in W\!L\backslash \pi_0(X^C)$
the number
$$
\sum_{\substack{W\!L_{y(C)}\cdot \sigma \in \\W\!L_{y(C)}\backslash \mor(y(C),x)}}
\frac{|\Gen(L)|}{|L| \cdot |(WC_{y(C)})_{\sigma}|},
$$
where $y(C)\co G/L \to X$ is some object in $\Pi_0(G,X)$ with
$ X^L(y) = C$ in $\pi_0(X^L)$. 

Recall the bijection $\alpha_F(x)$ from
\eqref{definition of alpha_F(x)}. In the case $F = \qq$ it becomes the map
\begin{multline*}
\alpha_{\qq}(x)\co \sum_{\substack{(L)\in \consub(G)\\ L \text{ finite cyclic}}} ~
\coprod_{W\!L\backslash \pi_0(X^L)} ~ W\!L_{y(C)}\backslash \mor(y(C),x) ~
\\
\xrightarrow{\cong} \con_{\qq}(H) = \{(K) \in \consub(H) \mid L \text{ cyclic}\}
\end{multline*}
which sends $W\!L\cdot \sigma \in W\!L\backslash \pi_0(X^L)$ to
$(\sigma(1L)^{-1} L\sigma(1L))$. 
Let
$$u_{[x]} \in \class_{\qq}(H)$$ 
be the element which assigns to
$(K) \in \con_{\qq}(H)$ the number
$\frac{|\Gen(L)|}{|L| \cdot |(WC_{y(C)})_{\sigma}|}$ if
$\sigma \in \mor(y(C),x)$ represents the preimage of $(K)$ under the 
bijection $\alpha_{\qq}(x)$. We conclude that
$\left(\alpha_{\qq}^G(X) \circ D^G(X) \circ \pr \circ \ch^G(X)\right)([x])$ is given by
the image under the structure map associated to the object $x\co G/H \to X$
$$\class_{\qq}(H) \to  \colim_{\Pi_0(G,X)} Q^*\underline{\class_{\qq}}$$
of the element $u_{[x]} \in \class_{\qq}(H)$ above.

Consider $(K) \in \con_{\qq}(H)$. Let $\sigma \in \mor(y(C),x)$ 
represent the preimage of $(K)$ under the bijection $\alpha_{\qq}(x)$.
Choose $g'$ such that $\sigma\co G/\langle g \rangle \to G/H$ is given by
$g''\langle g \rangle \mapsto g''g'H$.
Let $N_HK$ be the  normalizer in $H$ and $W_HK:= N_HK/K$ be the Weyl group of $K\subseteq H$.
Define a bijection 
$$f\co (W\!L_{y(C)})_{\sigma} \xrightarrow{\cong} W_H((g')^{-1}Lg'),
\hspace{5mm} g''L \mapsto (g')^{-1} g'' g'\cdot (g')^{-1}Lg'.$$ 
The map is well-defined because of
$$(W\!L_{y(C)})_{\sigma} = \{g''L  \in W\!L_{y(C)} \mid  (g')^{-1} g'' g' \in H\}$$
and the following calculation
\begin{multline*}
\left((g')^{-1} g'' g'\right)^{-1} (g')^{-1}Lg' (g')^{-1} g'' g'  ~ = ~
(g')^{-1} (g'')^{-1} g' (g')^{-1}Lg' (g')^{-1} g'' g'
\\
 ~ = ~ 
(g')^{-1} (g'')^{-1}L  g'' g' ~ = ~  (g')^{-1}Lg'.
\end{multline*}
One easily checks that $f$ is injective. Consider $h\cdot (g')^{-1}Lg'$
in $W_H((g')^{-1}Lg')$. Define $g_0 = g'h(g')^{-1}$. We have
$g_0 \in W\!L$. Since $h \cdot x(1H) = x(1H)$, we get $g_0L \in W\!L_y$.
Hence $g_0L$ is a preimage of $h\cdot (g')^{-1}Lg'$ under $f$.
Hence $f$ is bijective. This shows
$$|(W\!L_{y(C)})_{\sigma}| ~ = ~ |W_H((g')^{-1}Lg')|.$$
We conclude that $u_{[x]}$ is the element 
$$\con_{\qq}(H) = \{(K) \in \consub(H) \mid K \text{ finite cyclic}\} ~ \to \qq,
\hspace{5mm} (K) ~ \mapsto ~ \frac{|\Gen(K)|}{|K| \cdot |W_HK|}.$$
Since right multiplication with $|H|^{-1} \cdot \sum_{h \in H} h$ induces
an idempotent $\qq H$-linear map $\qq H \to \qq H$ whose image is $\qq$ with the trivial
$H$-action, the element $\HS_{\qq,H}([\qq]) \in \class_{\qq}(H)$ is given by
$$\con_{\qq}(H) ~ \to \qq,
\hspace{5mm} (K) ~ \mapsto ~ \frac{1}{|H|} \cdot |\{h \in H \mid \langle h \rangle \in (K)\}|.$$
From
\begin{eqnarray*}
|\{h \in H \mid \langle h \rangle \in (K)\}| & = & 
\left|\coprod_{K' \subseteq H,  K' \in  (K)} \Gen(K') \right|
\\
& = & \left|\{K' \subseteq H, K' \in  (K)\}\right| \cdot |\Gen(K)|
\\
& = & \frac{|H|}{|N_HK|} \cdot |\Gen(K)|
\\
& = &  \frac{|H| \cdot |\Gen(K)|}{|K| \cdot |W_HK|}
\end{eqnarray*}
we conclude
$$u_{[x]} ~ = ~ \HS_{\qq,H}([\qq]) \hspace{3mm} \in \class_{\qq}(H).$$
Now \eqref{the: commutative diagram for e^G_1(M): 1} and hence
Theorem \ref{the: commutative diagram for e^G_1(M)} follow. \qed


\tit{Examples}{Examples}

In this section we discuss some examples.
Recall that we have described the non-equivariant case in the introduction.

 
\subtit{Finite groups and connected non-empty fixed point sets}
{Finite groups and connected non-empty fixed point sets}

Next we consider the case where $G$ is a finite group, $M$ is  a
closed compact $G$-manifold, and 
$M^H$ is connected and non-empty for all subgroups $H \subseteq G$. 
Let $i \co \{\ast\} \to M$ be the $G$-map given by the inclusion of a
point into $M^G$. 
Since we have assumed that $M^G$ is connected, $i$ is unique up to
$G$-homotopy. Let $A(G)$ be the Burnside ring of formal differences
of finite $G$-sets (Definition \ref{def: Burnside ring}). 
We have the following commutative diagram:
$$\begin{CD}
U^G(M) @< U^G(i) < \cong < U^G(\{\ast\}) @< f_1 < \cong < A(G)
\\
@V e^G_1(M) VV @V e^G_1(\{\ast\}) VV @V j_1 VV 
\\
H_0^{\Or(G)}(M;\underline{R_{\qq}})
@<H_0^{\Or(G)}(i;\underline{R_{\qq}}) < \cong < 
H_0^{\Or(G)}(\{\ast\};\underline{R_{\qq}}) 
@<f_2 < \cong <
R_{\qq}(G)
\\
@V e^G_2(M) VV @V e^G_2(\{\ast\}) VV @V j_2 VV 
\\
H_0^{\Or(G)}(M;\underline{R_{\rr}})
@<H_0^{\Or(G)}(i;\underline{R_{\rr}}) < \cong < 
H_0^{\Or(G)}(\{\ast\};\underline{R_{\rr}}) 
@<f_3 < \cong <
R_{\rr}(G)
\\
@V e^G_3(M) VV @V e^G_3(\{\ast\}) V\cong V @V j_3 V\cong V 
\\
KO_0^G(M) @ <KO_0^G(i) << KO_0^G(\{\ast\}) @ < \id < \cong <
\,KO_0^G(\{\ast\}), 
\end{CD}
$$
where $j_1$ sends the class of a $G$-set $S$ in the Burnside ring $A(G)$ 
to the class of the rational $G$-representation
$\qq[S]$ and $j_2$ is the change-of-coefficients  homomorphism. The
homomorphism 
$KO_0^G(i)\co KO_0^G(\{\ast\}) \to KO_0^G(M)$ is split injective,
with a splitting given by the map $KO_0^G(\pr)$ induced by
$\pr\co M \to \{\ast\}$. The map $e^G_3(\{\ast\})$ is bijective
since the category $\Pi_0(G,\{\ast\})$ has $G \to \{\ast\}$ as
terminal object. This implies that 
$$e^G_3(M) \co H_0^{\Or(G)}(M;\underline{R_{\rr}}) \to KO_0^G(M)$$
is split injective. We have already explained in the Section
\ref{sec: The transformation eG(M)} 
that the map $j_2\co R_{\qq}(G) \to R_{\rr}(G)$ is rationally injective.
Since $R_{\qq}(G)$ is a torsion-free finitely generated abelian group,
$j_2$ is injective. Hence 
$$e^G_2(M) \co H_0^{\Or(G)}(M;\underline{R_{\qq}}) \to
H_0^{\Or(G)}(M;\underline{R_{\rr}})$$
is injective.  The upshot of this discussion is,
that $e^G_1(\chi^G(M))$ carries (integrally) the same information as $\Eul^G(M)$
because it is sent to $\Eul^G(M)$ by the injective map $e^G_3(M) \circ e^G_2(M)$. 

Analyzing the difference between  $e^G_1(\chi^G(M))$ and $\chi^G(M)$ is equivalent to analyzing
the map $j_1 \co A(G) \to \Rep_{\qq}(G)$, which sends $\chi^G(M)$ to the element given by
$\sum_{p \ge 0} (-1)^p \cdot [H_p(M;\qq)]$. Recall that $\chi^G(M) \in A(G)$ is given by
$$ \chi^G(M) ~ = ~ 
\sum_{(H) \in \consub(G)} \chi(W\!H\backslash M^H,W\!H\backslash M^{>H}) \cdot G/H]
~ = ~ \sum_{p \ge 0} (-1)^p \cdot \sharp_p(G/H),$$
where  $\chi(W\!H\backslash M^H,W\!H\backslash M^{>H}))$ is the non-equivariant
Euler characteristic and  $\sharp_p(G/H)$ is the number of equivariant cells
of the type $G/H \times D^p$ appearing in some $G$-$CW$-complex structure on $M$.
The following diagram commutes (see Theorem \ref{the: commutative diagram for e^G_1(M)})
\comsquare{A(G)}{j_1}
{R_{\qq}(G)}{\ch^G}{\HS_{\qq,G}}
{\prod_{(H) \in \consub(H)} \qq}{\pr}{\prod_{\substack{(H)
\consub(G)\\ H \text{ cyclic}}} \qq} 
where $\pr$ is the obvious projection, 
$\ch^G(S)$ has as entry for $(H) \subseteq \consub(G)$ the number
$\frac{1}{|W\!H|} \cdot |S^H|$ and $\HS_{\qq,G}(V)$ has as entry at $(H) \in \consub(G)$
for cyclic $H \subseteq G$ the number $\frac{\tr_{\qq}(l_h)}{|W\!H|}$, where
$\tr(l_h) \in \qq$ is the trace of the endomorphism of the rational 
vector space $V$ given by multiplication with $h$ for some generator
$h \in H$. 
The vertical arrows $\ch^G$ and $\HS_{\qq,G}$ are rationally bijective
and $\chi^G(M)(\chi^G(M))$ has as component belonging to $(H)\in
\consub(G)$ the number 
$\chi^{\qq W\!H}(M^H) = \frac{\chi(M^H)}{|W\!H|}$
(see Lemma \ref{lem: injectivity of the character map chi^G} and Lemma
\ref{R_F(H) and class_F(H)}). This implies that $\ch^G$ and $\HS_{\qq,G}$ 
are injective because their sources are torsion-free finitely
generated abelian groups. 
Moreover, $\chi^G(M) \in A(G)$ carries integrally the same information as all
the collection of  Euler characteristics $\{\chi(M^H) \mid (H) \in
\consub(H)\}$, 
whereas $j_1(\chi^G(M)) = \sum_{p \ge 0} (-1)^p \cdot [H_p(M;\qq)]$
carries integrally the same information as  the collection of the
Euler characteristics 
$\{\chi(M^H) \mid (H) \in \consub(H), H\text{ cyclic}\}$. In particular
$j_1$ is injective if and only if $G$ is cyclic. But any element 
$u$ in $A(G)$ occurs as $\chi^G(M)$ for a closed smooth $G$-manifold $M$ 
for which $M^H$ is connected and non-empty for all $H \subseteq G$
(see \cite[Section 7]{Lueck-Rosenberg(2002a)}). Hence, given a finite
group $G$, the elements $\chi^G(M)$ and $j_1(\chi^G(M))$ carry the same
information for all such $G$-manifolds $M$ if and only if $G$ is cyclic.

From this discussion we conclude 
that $\Eul^G(M)$ does not carry torsion information 
in the case where $G$ is finite and $M^H$ 
is connected and non-empty for all $H \subseteq G$, since
$R_{\rr}(G)$ is a torsion-free finitely generated abelian group.
This is different from the case where one allows non-connected fixed
point sets, as the following example shows.

 
\subtit{The equivariant Euler class carries torsion information}
{EulG(M) carries torsion information}

Let $S_3$ be the symmetric group on $3$ letters. It has the presentation
$$S_3 = \langle s,t \mid s^2 = 1, t^3 = 1, sts= t^{-1}\rangle.$$ 
Let $\rr$ be the trivial $1$-dimensional real representation of $S_3$.
Denote by $\rr^-$ the one-dimensional real representation on which $t$
acts trivially and $s$ acts by $-\id$. Denote by $V$ the
$2$-dimensional irreducible real representation of $S_3$; we can
take $\rr^2$ as the underlying
real vector space of $V$, with $s \cdot (r_1,r_2) = (r_2,r_1)$ and
$t \cdot (r_1,r_2) = (-r_2,r_1 - r_2)$. Then $\rr$, $\rr^-$ and $V$ are the 
irreducible real representations of $S_3$, and 
$\rr S_3$ is as an $\rr S_3$-module isomorphic to $\rr \oplus \rr^-
\oplus V \oplus V$. Let $L_2$ be the cyclic group of order two
generated by $s$ and let
$L_3$ be the cyclic group of order three  generated by $t$. Any finite
subgroup of 
$S_3$ is conjugate to precisely one of the subgroups $L_1=\{1\}$, $L_2$,
$L_3$ or $S_3$. One easily checks that
 $V^{L_2} \cong \rr$, $V^{L_3} \cong 0$, $(\rr^-)^{L_2} = 0$ and
$(\rr^-)^{L_3} \cong \rr$ as real vector spaces. 
Put $W = \rr \oplus \rr^- \oplus V$. Then $W^{S_3}  \cong \rr$,
$W^{L_2} \cong W^{L_3} \cong \rr^2$ and $W \cong \rr^4$ as real vector spaces.
Let $M$ be the closed $3$-dimensional $S_3$-manifold $SW$. Then
\begin{eqnarray*}
M & \cong & S^3;
\\
M^{L_2} & \cong & S^1;
\\
M^{L_3} & \cong & S^1;
\\
M^{S_3} & \cong & S^0.
\end{eqnarray*}

Since $\chi(M^{S_3})\ne 0$, $\chi^G(M)\in U^G(M)$ cannot vanish.
But since the fixed sets for all cyclic subgroups have vanishing Euler
characteristic, Theorem \ref{the: rational information carried by Eul^G(M)}
implies that $\Eul^G(M)$ is a torsion element 
in $KO_0^G(M)$. We want to show that it has order precisely two.

Let $x_i\co S_3/L_i \to X$ for $i = 1,2,3$ be a $G$-map. Let $x_-\co
S_3/S_3 \to X$ and 
$x_+\co S_3/S_3 \to X$ be the two different $G$-maps for which
$M^{S_3}$ is the union of the images of $x_-$ and $x_+$. Then
$x_1$, $x_2$, $x_3$, $x_-$ and $x_+$ form a complete set of representatives
for the isomorphism classes of objects in $\Pi_0(S_3,M)$.  Notice for $i = 1,2,3$ that 
$\mor(x_i,x_-)$ and $\mor(x_i,x_+)$ consist of precisely one element. 
Therefore we get an exact sequence
\begin{multline}
R_{\rr}(\zz/2) \oplus R_{\rr}(\zz/3)
\xrightarrow{\squarematrix{i_2}{i_3}{-i_2}{-i_3}} 
R_{\rr}(G) \oplus R_{\rr}(G) 
\\
\xrightarrow{s_- + s_+} 
\colim_{\Pi_0(S_3,M)} Q^*\underline{R_{\rr}} \to 0,
\label{exact sequence for colim for S_3-manifold M}
\end{multline}
where $i_2\co R_{\rr}(\zz/2) \to R_{\rr}(S_3)$ and
$i_3\co R_{\rr}(\zz/3) \to R_{\rr}(S_3)$ are  the induction homomorphisms associated to
any injective group homomorphism from $\zz/2$ and $\zz/3$ into $S_3$ and
$s_-$ and $s_+$ are the structure maps of the colimit belonging to the
objects $x_-$ and $x_+$. Define a map
$$\delta\co R_{\rr}(G) \to \zz/2, \hspace{5mm} 
\lambda_{\rr} \cdot [\rr] + \lambda_{\rr^-} \cdot [\rr^-] +
\lambda_{V} \cdot [V]  
~ \mapsto \overline{\lambda_{\rr}}   +  \overline{\lambda_{\rr^-}} + 
\overline{\lambda_V}.$$
If $\rr$ denotes the trivial and $\rr^-$ denotes the non-trivial
one-dimensional real $\zz/2$-representation, then 
\[
\begin{aligned}
i_2&\co R_{\rr}(\zz/2) \to R_{\rr}(S_3), \\
&\mu_{\rr} \cdot [\rr] + \mu_{\rr^-} \cdot [\rr^-] ~ \mapsto 
\mu_{\rr} \cdot [\rr] +  \mu_{\rr^-} \cdot [\rr^-] + 
(\mu_{\rr} +  \mu_{\rr^-}) \cdot [V].
\end{aligned}
\]
If $\rr$ denotes the trivial one-dimensional and $W$ the
$2$-dimensional irreducible real $\zz/3$-representation, then 
$$i_3\co R_{\rr}(\zz/3) \to R_{\rr}(S_3) \hspace{5mm}
\mu_{\rr} \cdot [\rr] + \mu_{W} \cdot [W] ~ \mapsto 
\mu_{\rr} \cdot [\rr] +  \mu_{\rr} \cdot [\rr^-] + 
2 \mu_{W} \cdot [V].$$
This implies that the following sequence is exact
\begin{eqnarray}
&R_{\rr}(\zz/2) \oplus R_{\rr}(\zz/3) \xrightarrow{i_2 + i_3}
R_{\rr}(G) \xrightarrow{\delta} \zz/2\to 0. &
\label{exact sequence for the Artin defect of R_{rr}(S_3)}
\end{eqnarray}
We conclude from the exact sequences 
\eqref{exact sequence for colim for S_3-manifold M} and
\eqref{exact sequence for the Artin defect of R_{rr}(S_3)} above 
that the epimorphism
$$s_- + s_+ \co  R_{\rr}(G) \oplus R_{\rr}(G) \to
\colim_{\Pi_0(S_3,M)} Q^* \underline{R_{\rr}}$$
factorizes through the map
$$u:= \squarematrix{1}{1}{0}{\delta} \co  
R_{\rr}(G) \oplus R_{\rr}(G) \to R_{\rr}(G) \oplus \zz/2$$
to an isomorphism
\begin{eqnarray}
 v \co R_{\rr}(G) \oplus \zz/2 & \xrightarrow{\cong} &
\colim_{\Pi_0(S_3,M)} Q^*\underline{R_{\rr}}.
\label{computation of colim_{z colon S_3/H to M  in Pi_0(S_3,M)} R_{rr}(H)}
\end{eqnarray}
Define a map
$$f\co U^{S_3}(M) \to R_{\rr}(G) \oplus \zz/2$$
by 
\begin{eqnarray*}
f([x_1]) & := & = ([\rr[S_3]],0);
\\
f([x_2]) & := & = ([\rr[S_3/L_2]],0);
\\
f([x_3]) & := & = ([\rr[S_3/L_3]],0);
\\
f([x_-]) & := & = ([\rr],0);
\\
f([x_+]) & := & = ([\rr],1).
\end{eqnarray*}
The reader may wonder why $f$ does not look symmetric in $x_-$ and $x_+$. This comes from the 
choice of $u$ which affects the isomorphism $v$. The composition
$$U^G(M) \xrightarrow{e^G_1(M)} H_0^{\Or(G)}(M;\underline{R_{\qq}})
\xrightarrow{e^G_2(M)} H_0^{\Or(G)}(M;\underline{R_{\rr}})$$
agrees with the composition
\begin{multline*}U^G(M) \xrightarrow{f}  R_{\rr}(G) \oplus \zz/2 
\xrightarrow{v}  
\colim_{\Pi_0(S_3,M)} Q^*\underline{R_{\rr}}
\xrightarrow{\beta^G_{\rr}(M)} 
H_0^{\Or(G)}(M;\underline{R_{\rr}}).
\end{multline*}

We get
\begin{eqnarray}
\chi^G(M) & = & [x_+] + [x_-]  - 2 \cdot [x_2] - [x_3]  + [x_1]
\hspace{5mm} \in U^G(M)
\label{computation of chi^G(M)}
\end{eqnarray} 
since the image of the element on the right side and the image of
$\chi^G(M)$ under the injective character map $\ch^G(M)$ (see 
\eqref{definition of character map chi^G(X)} 
and Lemma \ref{lem: injectivity of the character map chi^G}) agree 
by the following calculation
\begin{eqnarray*}
\ch^G(\chi^G(M))   & = & 
 \frac{\chi(M^{S_3}(x_-))}{|(W\!S_3)_{x_-}|}  \cdot [x_-] + 
 \frac{\chi(M^{S_3}(x_+))}{|(W\!S_3)_{x_+}|}  \cdot [x_+] + 
\frac{\chi(M^{L_2})}{|W\!L_2|}  \cdot [x_2] 
\\ & & \hspace{20mm} + \frac{ \chi(M^{L_3})}{W\!L_3} \cdot [x_3] +
 \frac{ \chi(M)}{W\!L_1} \cdot [x_1]
\\
& = &
[x_-] +  [x_+]
\\
& = & \ch^G(M)( [x_+] + [x_-]  - 2 \cdot [x_2] - [x_3]  + [x_1]).
\end{eqnarray*}
Now one easily checks
\begin{eqnarray}
f(\chi^G(M)) & = & (0,1) \hspace{5mm} \in R_{\rr}(S_3) \oplus \zz/2.
\label{f(chi^G(M))}
\end{eqnarray}
Since $R_{\rr}(G) \oplus \zz/2  \xrightarrow{v}  
\colim_{\Pi_0(S_3,M)} Q^* \underline{R_{\rr}}
\xrightarrow{\beta^G_{\rr}(M)} 
H_0^{\Or(G)}(M;\underline{R_{\rr}})$ is a composition of isomorphisms,
we conclude 
that $e^G_2(M) \circ e^G_1(M)(\chi^G(M))$ is an element of order two in
$H_0^{\Or(G)}(M;\underline{R_{\rr}})$. Since $\dim(M)  \le 4$
and $\dim(M^{>1})\le 2$, we
conclude from 
Theorem \ref{the: e^G(X)}  \eqref{the: e^G(X): injectivity of e^G_3(X)
for small dim} 
that $e^G_3(M)\co H_0^{\Or(G)}(M;\underline{R_{\rr}}) \to KO_0^G(M)$
is injective. 
We conclude from Theorem \ref{the: e^G(M)(chi^G(M)) = Eul^G(M)} and 
Theorem \ref{the: e^G(X)}
\eqref{the: e^G(X): e^G(X) = e^G_3(X) circ e_2^G(X) circ e_1^G(X)}
that $\Eul^G(M) \in KO_0^G(M)$ is an element of order two as promised in
Theorem \ref{the: example of Eul^G(M) of order two}.

 
\subtit{The equivariant Euler class is independent of the stable
equivariant Euler characteristic}{Independence of stable
equivariant Euler characteristic}

In this subsection we will give examples to show that $\Eul^G(M)$ 
is independent of the stable equivariant Euler characteristic with
values in the Burnside ring $A(G)$, in the sense that it is possible for
either one of these invariants to vanish while the other does not
vanish.

For the first example, take $G=\bZ/p$ cyclic of prime order,
so that $G$ has only two subgroups (the trivial subgroup and $G$
itself) and $A(G)$ has rank $2$.
We will see that it is possible for $\Eul^G(M)$ to be
non-zero, even rationally, while $\chi^G_s(M)=0$ in $A(G)$ (see
Definition \ref{def: more Euler characteristics}).
To see this, we will construct a closed $4$-dimensional $G$-manifold
$M$ with $\chi (M)=0$ and such that 
$M^G$ has two components of dimension $2$,
one of which is $S^2$ and the other of which is a surface $N^2$ of
genus $2$, so that $\chi(N)=-2$. Then 
\[
\chi(M^G)=\chi(S^2\amalg N^2) = \chi(S^2) + \chi(N^2) = 2-2 =0,
\]
while also $\chi(M^{\{1\}})=
\chi(M)=0$, so that $\chi^G_s(M)=0$ in $A(G)$ and hence also $\chi^G_a(M)=0$.

For the construction, simply choose any bordism $W^3$ between $S^2$
and $N^2$, and let
\[
M^4 = \left(S^2\times D^2\right) \cup_{S^2\times S^1}
\left(W^3\times S^1\right) \cup_{N^2\times S^1}
\left(N^2\times D^2\right) .
\]

We give this the $G$-action which is trivial on the $S^2$, $W^3$, and
$N^2$ factors, and which is rotation by ${2\pi}/{p}$ on
the $D^2$ and $S^1$ factors. Then the action of $G$ is free except
for $M^G$, which consists of $S^2\times \{0\}$ and of  $N^2\times \{0\}$.
Furthermore, we have
\begin{eqnarray*}
\chi(M) &=& \chi \left(S^2\times D^2\right) -\chi \left(S^2\times S^1\right)
+\chi\left(W^3\times S^1\right)\\
&&\qquad -\chi \left(N^2\times S^1\right)
+\chi \left(N^2\times D^2\right) \\
&=&2-0+0-0-2=2-2=0.
\end{eqnarray*}
Thus $\chi(M)=\chi(M^G)=0$ and $\chi_s^G(M)=0$ in $A(G)$. On the other
hand, $\Eul^G(M)$ is non-zero, even rationally, since from it 
(by Theorem \ref{the: rational information carried by Eul^G(M)}) we can
recover the two (non-zero) Euler characteristics of the two components
of $M^G$. 

For the second example, take $G=S_3$ and retain the
notation of Subsection \ref{subsec: EulG(M)
carries torsion information}. By \cite[Theorem
7.6]{Lueck-Rosenberg(2002a)}, there is a closed $G$-manifold 
$M$ with $M^H$ connected for each subgroup $H$ of $G$, with
$\chi(M^H)=0$ for $H$ cyclic, and with $\chi(M^G)\ne 0$. (Note
that $G$ is the only noncyclic subgroup of $G$.) In fact, we can write
down such an example explicitly; simply let $Q=W'\oplus \bR\oplus\bR$,
the $S_3$-representation $\bR\oplus\bR\oplus
\bR\oplus\bR^-\oplus V$, and let $M=SW'$ be the unit ball in $W'$.
Then each fixed set in $M$ is a sphere of dimension bigger by $2$ than
in the example of \ref{subsec: EulG(M)
carries torsion information}, so $M\cong S^5$, $M^{L_2}\cong S^3$,
$M^{L_3}\cong S^3$, and $M^{G}\cong S^2$. Since the fixed sets are all
connected and each fixed set of a cyclic subgroup has vanishing Euler
characteristic, it follows by Subsection \ref{subsec: Finite groups
and connected non-empty fixed point sets} that $\Eul^G(M)=0$. On the
other hand, since $\chi(M^{G})=2$, $\chi_s^G(M)\ne 0$ in $A(G)$.

 
\subtit{The image of  the equivariant Euler class under the assembly
maps}{The image of EulG(M) under the assembly maps} 

Now let us consider an infinite (discrete) group $G$. Let
$\underline{E}G$ be a 
model for the \emph{classifying space for proper $G$-actions}, i.e., a
$G$-$CW$-complex $\underline{E}G$ such that $\underline{E}G^H$ is
contractible (and in particular non-empty) for finite $H \subseteq G$
and $\underline{E}G^H$ is empty for infinite $H \subseteq G$. It has the
universal property that for any proper $G$-$CW$-complex $X$ there is up to
$G$-homotopy precisely one $G$-map $X \to \underline{E}G$. This
implies that all models for  $\underline{E}G$ are $G$-homotopy
equivalent. If $G$ is a word-hyperbolic group,
its Rips complex is a model for $\underline{E}G$
\cite{Meintrup(2000)}.
If $G$ is a discrete subgroup of the Lie group $L$ with finitely many
path components, 
then $L/K$ for a maximal compact subgroup $K$ with the (left)
$G$-action is a model for 
$\underline{E}G$ (see \cite[Corollary 4.14]{Abels (1978)}).
If $G$ is finite, $\{\ast\}$ is a model for $\underline{E}G$.

Consider a proper smooth $G$-manifold $M$. Let $f\co M \to \underline{E}G$
be a $G$-map. We obtain a commutative diagram
$$\begin{CD}
U^G(M) @> U^G(f) >> U^G(\underline{E}G) @> \id >> U^G(\underline{E}G)
\\
@V e^G_1(M) VV @V e^G_1(\underline{E}G) VV @V j_1 VV 
\\
H_0^{\Or(G)}(M;\underline{R_{\qq}}) @>H_0^{\Or(G)}(f;\underline{R_{\qq}}) >>
H_0^{\Or(G)}(\underline{E}G;\underline{R_{\qq}}) 
@> \asmb_1 >> K_0(\qq G)
\\
@V e^G_2(M) VV @V e^G_2(\underline{E}G) VV @V j_2 VV 
\\
H_0^{\Or(G)}(M;\underline{R_{\rr}}) @>H_0^{\Or(G)}(f;\underline{R_{\rr}})>>
H_0^{\Or(G)}(\underline{E}G;\underline{R_{\rr}}) 
@> \asmb_2 >> K_0(\rr G)
\\
@V e^G_3(M) VV @V e^G_3(\underline{E}G) VV @V j_3 VV 
\\
KO_0^G(M) @ >KO_0^G(f) >> KO_0^G(\underline{E}G) @> \asmb_3 >>
KO_0(C_r^*(G;\rr)) 
\end{CD}
$$
Here $\asmb_i$ for $i = 1,2,3$ are the assembly maps appearing in the
Farrell-Jones 
Isomorphism Conjecture and the Baum-Connes Conjecture. The maps 
$\asmb_i$ for $i = 1,2$ are the obvious maps 
$\colim_{\Orcat} \underline{R_F} \to K_0(FG)$ for $F = \qq,\rr$
under the identifications $R_F(H) = K_0(FH)$ for finite $H \subseteq G$
and $\beta^G_F(M)$ of \eqref{identification beta^G_1(X) of colimit and
H_0^{Or(G)}}. 
The Baum-Connes assembly map is given by the index with values in the reduced
real group $C^*$-algebra $C_r^*(G;\rr)$. The Farrell-Jones
Isomorphism Conjecture and the Baum-Connes Conjecture predict that
$\asmb_i$ for $i = 1,2,3$ are bijective. The abelian group
$U^G(\underline{E}G)$ is the 
free abelian group with $\{(H) \in \consub(G) \mid |H| < \infty\}$ as
basis and the map 
$j_1$ sends the basis element $(H)$ to the class of the finitely
generated projective 
$\qq G$-module $\qq[G/H]$. The maps $j_2$ and $j_3$ are change-of-rings
homomorphisms. The maps $e^G_2(M)$ and $e^G_3(M)$ are rationally injective
(see Theorem \ref{the: e^G(X)} 
\eqref{the: e^G(X): rational injectivity of e^G_2(X)}
and \eqref{the: e^G(X): Chern character}). 
If $M^H$ is connected and non-empty for all finite subgroups $H\subseteq G$,
then the horizontal arrows $U^G(f)$, $H_0^{\Or(G)}(f;\underline{R_{\qq}})$ and
$H_0^{\Or(G)}(f;\underline{R_{\rr}})$ are bijective.

In contrast to the case where $G$
is finite, the groups $K_0(\qq G)$, $K_0(\rr G)$ and
$KO_0(C_r^*(G;\rr))$ may not be 
torsion free. The problem whether $\asmb_i$ is bijective for $i = 1,2$ or $3$ 
is a difficult and  in general unsolved problem. Moreover,
$\underline{E}G$ is a complicated $G$-$CW$-complex for
infinite $G$, whereas for a finite group $G$ we can take
$\{\ast\}$ as a model for $\underline{E}G$.


 \end{document}

%% file: gtspec.tex
\def\ifplaintex{\expandafter\ifx\csname documentclass\endcsname\relax}

\expandafter\ifx\csname epsfbox\endcsname\relax\input epsf\fi

\ifplaintex 
\else
\fi

\def\gt{{\mathsurround=0pt\it $\cal G\mskip-2mu$eometry \&\ 
$\cal T\!\!$opology}}        

\def\gtp{{\mathsurround=0pt\it $\cal G\mskip-2mu$eometry \&\ 
$\cal T\!\!$opology $\cal P\!$ublications}}  

\def\lognumber#1{\def\thelognumber{#1}}
\def\volumenumber#1{\def\thevolumenumber{#1}}
\def\papernumber#1{\def\thepapernumber{#1}}
\def\volumeyear#1{\def\thevolumeyear{#1}}

\def\pagenumbers#1#2{\def\startpage{#1}\def\finishpage{#2}}
\def\published#1{\def\publishdate{#1}}
\def\proposed#1{\def\theproposer{#1}}
\def\seconded#1{\def\theseconders{#1}}
\def\received#1{\def\receiveddate{#1}}

\def\accepted#1{\def\accepteddate{#1}}
\def\asciititle#1{\def\theasciititle{#1}}

\def\coverauthors#1{\def\thecoverauthors{#1}}
\def\asciiauthors#1{\def\theasciiauthors{#1}}
\def\asciiaddress#1{\def\theasciiaddress{#1}}
\def\asciiemail#1{\def\theasciiemail{#1}}
\long\def\asciiabstract#1{\long\def\theasciiabstract{#1}}
\def\asciikeywords#1{\def\theasciikeywords{#1}}

\let\\\par\let\thevolumenumber\relax\let\thepapernumber\relax
\let\thevolumeyear\relax\let\thesamplenumber\relax\let\startpage\relax
\let\finishpage\relax\let\publishdate\relax\let\receiveddate\relax
\let\reviseddate\relax\let\accepteddate\relax\let\theasciititle\relax
\let\theasciiauthors\relax\let\theasciiaddress\relax
\let\theasciiabstract\relax\let\theasciikeywords\relax
\let\theasciiemail\relax\let\theshortauthors\relax\let\theshorttitle\relax
\let\thecoverauthors\relax

\long\def\maketitlep{   

\count0=\startpage

\gt\hfill      
\hbox to 77pt{\vbox to 0pt{\vglue -15pt\epsfbox{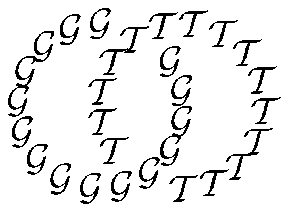}\vss}\hss}
\break
{\small\ifx\thesamplenumber\relax 
Volume \else Sample
\fi\thevolumenumber\ (\thevolumeyear)
\startpage--\finishpage\nl
Published: \publishdate}
\vglue 0.5truein plus 0.4fil minus 0.1truein

{\parskip=0pt\leftskip 0pt plus 1fil\def\\{\par\smallskip}{\ifplaintex\large
\else\Large\fi\bf\thetitle}\par\medskip}   

\vglue 0pt plus 0.1fil 

{\parskip=0pt\leftskip 0pt plus 1fil\def\\{\par}{\sc\theauthors}
\par\medskip}

\vglue 0pt plus 0.1fil 

{\small\parskip=0pt\let\newline\\
{\leftskip 0pt plus 1fil\def\\{\par}{\sl\theaddress}\par}
\expandafter\ifx\theemail\relax    
\relax\else\vglue 5pt plus 0.02fil minus 2pt\def\\{\stdspace{\rm 
and}\stdspace} 
\cl{Email:\stdspace\tt\theemail}\fi
\ifx\theurl\relax                  
\relax\else\vglue 5pt plus 0.02fil minus 2pt\def\\{\stdspace{\rm 
and}\stdspace}
\cl{URL:\stdspace\tt\theurl}\fi\par}

\vglue 7pt plus 0.3fil minus 3pt

{\bf Abstract}
\vglue 5pt plus 0.1fil minus 2pt

\theabstract

\vglue 7pt plus 0.3fil minus 3pt

{\bf AMS Classification numbers}\quad Primary:\quad \theprimaryclass

Secondary:\quad \thesecondaryclass

\vglue 5pt plus 0.3fil minus 2pt

{\bf Keywords:}\quad \thekeywords

\vglue 10pt plus 0.5fil minus 5pt

{\small  Proposed: \theproposer\hfill Received: \receiveddate\nl
Seconded: \theseconders\hfill 
\ifx\reviseddate\relax                         
Accepted: \accepteddate                        
\else
Revised: \reviseddate                          
\fi}
\eject
}       

\let\maketitlepage\maketitlep
\let\maketitle\maketitlepage

\font\phead=cmsl9 scaled 950
\font=cmsl9 scaled 1050
\font\pnum=cmbx10 scaled 913
\font\lnum=cmbx10 
\font\pfoot=cmsl9 scaled 950
\font=cmsl9 scaled 1050
\ifplaintex
\headline{\vbox to 0pt{\vskip -4.5mm\line{\small\phead\ifnum
\count0=\startpage ISSN 1364-0380 (on line)
1465-3060 (printed) \hfill {\pnum\folio}\else\ifodd\count0\def\\{ }%
\ifx\theshorttitle\relax\thetitle\else\theshorttitle\fi\hfill{\pnum\folio}
\else\def\\{ and }{\pnum\folio}\hfill\ifx\theshortauthors\relax\theauthors
\else\theshortauthors\fi\fi\fi}\vss}}
\footline{\vbox to 0pt{\vglue 0mm\line{\small\pfoot\ifnum\count0=\startpage
\copyright\ \gtp\hfill\else
\gt, Volume \thevolumenumber\ (\thevolumeyear)\hfill\fi}\vss
}}
\else
\makeatletter
\def\@oddhead{{\small\lhead\ifnum\count0=\startpage ISSN 1364-0380 (on line)
1465-3060 (printed) \hfill {\lnum\number\count0}\else\ifodd\count0
\def\\{ }\ifx\theshorttitle\relax \thetitle \else\theshorttitle\fi\hfill
{\lnum\number\count0}\else\def\\{ and }{\lnum\number\count0}
\hfill\ifx\theshortauthors\relax 
\theauthors\else\theshortauthors\fi\fi\fi}}\def\@evenhead{@oddhead}
\def\@oddfoot{\small\lfoot\ifnum\count0=\startpage\copyright\ \gtp\hfill\else
\gt, Volume \thevolumenumber\ (\thevolumeyear)\hfill\fi}
\def\@evenfoot{@oddfoot}
\makeatother
\fi

\newwrite\gtoutfile
\long\gdef\makeheadfile{  
{\def\\{, }\def\s{ }
\immediate\openout\gtoutfile head.xxx
\immediate\write\gtoutfile{To: math@arxiv.org}
\immediate\write\gtoutfile{Subject: put OR rep NNNNN:pppp}
\immediate\write\gtoutfile{--text follows this line--}
\immediate\write\gtoutfile{Proxy-for: \ifx\theasciiauthors\relax
\theauthors\else\theasciiauthors\fi\s<\ifx\theasciiemail\relax\theemail\else\theasciiemail\fi>}
\immediate\write\gtoutfile{\noexpand\\}
\immediate\write\gtoutfile{Authors: \ifx\theasciiauthors\relax
\theauthors\else\theasciiauthors\fi}
{\def\\{ }\immediate\write\gtoutfile{Title: \ifx\theasciititle\relax
\thetitle\else\theasciititle\fi}}
\immediate\write\gtoutfile{Subj-class: GT or GR or SG or ...}
\immediate\write\gtoutfile{MSC-class: \theprimaryclass\ifx\thesecondaryclass\relax\else, \thesecondaryclass\fi}
\immediate\write\gtoutfile{Journal-ref: Geom. Topol. \thevolumenumber\s
(\thevolumeyear) \startpage-\finishpage}
\immediate\write\gtoutfile{Comments: Published in Geometry and Topology at}
\immediate\write\gtoutfile{    http://www.maths.warwick.ac.uk/gt/GTVol\thevolumenumber/paper\thepapernumber.abs.html}
\immediate\write\gtoutfile{\noexpand\\}
\immediate\write\gtoutfile{}
\ifx\theasciiabstract\relax
\immediate\write\gtoutfile{\theabstract}\else
\immediate\write\gtoutfile{\theasciiabstract}\fi
\immediate\write\gtoutfile{}
\immediate\write\gtoutfile{\noexpand\\}
\immediate\write\gtoutfile{}
\immediate\closeout\gtoutfile}}  

\def\maketitlepage{\maketitlep\makeheadfile}
\let\maketitle\maketitlepage